\documentclass{gtmon_a}
\pdfoutput=1

\proceedingstitle{The Zieschang Gedenkschrift}
\conferencestart{5 September 2007}
\conferenceend{8 September 2007}
\conferencename{Conference in honour of Heiner Zieschang}
\conferencelocation{Toulouse, France}

\editor{Michel Boileau}
\givenname{Michel}
\surname{Boileau}

\editor{Martin Scharlemann}
\givenname{Martin}
\surname{Scharlemann}

\editor{Richard Weidmann}
\givenname{Richard}
\surname{Weidmann}

\title[Remarks on the cohomology of finite fundamental groups of $3$--manifolds]{Remarks on the cohomology\\ of finite fundamental groups of $3$--manifolds}  

\author{Satoshi Tomoda}
\givenname{Satoshi}
\surname{Tomoda}
\address{
Department of Mathematics and Statistics\\
Okanagan College\\\newline
1000 KLO Road\\
Kelowna, B.C.
V1Y 4X8\\Canada}
\email{STomoda@Okanagan.bc.ca}
\urladdr{}

\author{Peter Zvengrowski}
\givenname{Peter}
\surname{Zvengrowski}
\address{
Department of Mathematics and Statistics\\
University of Calgary\\\newline
Calgary T2N 1N4\\Canada
}
\email{zvengrow@ucalgary.ca}
\urladdr{}

\volumenumber{14}
\issuenumber{}
\publicationyear{2008}
\papernumber{023}
\startpage{519}
\endpage{556}

\doi{}
\MR{}
\Zbl{}

\arxivreference{} 


\received{6 August 2006}
\revised{}
\accepted{14 February 2007}
\published{29 April 2008}
\publishedonline{29 April 2008}
\proposed{}
\seconded{}
\corresponding{}
\version{}

\makeatletter
\def\cnewtheorem#1[#2]#3{\newtheorem{#1}{#3}[section]
\expandafter\let\csname c@#1\endcsname\c@Thm}
\makeatother

\let\xysavmatrix\xymatrix
\def\xymatrix{\disablesubscriptcorrection\xysavmatrix}
\AtBeginDocument{\let\bar\wbar\let\tilde\wtilde\let\hat\what}
\AtBeginDocument{\newcommand{\<}{\langle}
\renewcommand{\>}{\rangle}}
\def\S{Section }
\makeop{hom}
\makeop{Spin}
\makeop{cat}

\newtheorem{Thm}{Theorem}[section]
\cnewtheorem{Cor}[Thm]{Corollary}
\cnewtheorem{Lem}[Thm]{Lemma}
\cnewtheorem{Pro}[Thm]{Proposition}

\theoremstyle{definition}
                             \newtheorem{Rem}[Thm]{Remark}
                             
\cnewtheorem{Def}[Thm]{Definition}

\makeautorefname{Thm}{Theorem}
\makeautorefname{Cor}{Corollary}
\makeautorefname{Pro}{Proposition}

\newcommand{\ZZ}{\mathbb {Z}}
\newcommand{\RR}{\mathbb {R}}

\newcommand{\Proof} {\begin{proof}}

\newcommand\nc{\newcommand}
\nc\surj{\twoheadrightarrow}
\nc\onto{\surj}
\nc\inj{\rightarrowtail}
\nc\iso{\approx}
\nc\inc{\hookrightarrow}
\nc\homeo{\cong}
\nc\lto{\longrightarrow}
\nc\ve{\varepsilon}
\nc\ctimesc{{\mathcal C}\otimes{\mathcal C}}
\nc{\bm}[1]{{\boldmath {$#1$}}}
\nc{\st}{\rule[-2pt]{0pt}{12pt}}

\begin{document}

\begin{asciiabstract}
Computations based on explicit 4-periodic resolutions are given for
the cohomology of the finite groups G known to act freely on S^3,
as well as the cohomology rings of the associated 3-manifolds
(spherical space forms) M = S^3/G.  Chain approximations to the
diagonal are constructed, and explicit contracting homotopies also
constructed for the cases G is a generalized quaternion group, the
binary tetrahedral group, or the binary octahedral group.  Some
applications are briefly discussed.  
\end{asciiabstract}

\begin{htmlabstract}
Computations based on explicit 4&ndash;periodic resolutions are given for
the cohomology of the finite groups G known to act freely on S<sup>3</sup>,
as well as the cohomology rings of the associated 3&ndash;manifolds
(spherical space forms) M = S<sup>3</sup>/G.  Chain approximations to the
diagonal are constructed, and explicit contracting homotopies also
constructed for the cases G is a generalized quaternion group, the
binary tetrahedral group, or the binary octahedral group.  Some
applications are briefly discussed.
\end{htmlabstract}

\begin{abstract}
Computations based on explicit $4$--periodic resolutions are given for
the cohomology of the finite groups $G$ known to act freely on $S^3$,
as well as the cohomology rings of the associated $3$--manifolds
(spherical space forms) $M = S^3/G$.  Chain approximations to the
diagonal are constructed, and explicit contracting homotopies also
constructed for the cases $G$ is a generalized quaternion group, the
binary tetrahedral group, or the binary octahedral group.  Some
applications are briefly discussed.  
\end{abstract}

\maketitle

\allowdisplaybreaks 

\section{Introduction} \label{sec:1}
The structure of the cohomology rings of $3$--manifolds is an area to which
Heiner Zieschang devoted much work and energy, especially from $1993$
onwards.  This could be considered as part of a larger area of his interest,
the degrees of maps between oriented $3$--manifolds, especially
the existence of degree one
maps, which in turn have applications in unexpected areas such as relativity
theory (cf Shastri, Williams and Zvengrowski \cite{swz} and Shastri and Zvengrowski \cite{sz}).  References 
\cite{adhsz,bhzz1,bhzz2,hkwz,hwz1,hwz2,hwz3,hz2,hz1} in this paper, all
involving work of Zieschang, his students Aaslepp, 
Drawe, Sczesny, and various
colleagues, attest to his enthusiasm for these topics and the remarkable
energy he expended studying them.

Much of this work involved Seifert manifolds, in particular, references
\cite{adhsz,bhzz1,bhzz2,hkwz,hwz2,hz1}.  
Of these, \cite{bhzz1,bhzz2,hz1} (together with \cite{blpz,bz})
successfully completed the programme of computing the ring structure
$H^*(M)$ for any orientable Seifert manifold $M$ with $G := \pi_1(M)$
infinite.  
Any such Seifert manifold $M$ (apart from $S^1 \times S^2$ and
$\RR P^3 \# \RR P^3$) is irreducible, hence aspherical (ie, an
Eilenberg--MacLane space $K(G, 1)$) by a well known application of the
Papakyriakopolous sphere theorem (see Hempel \cite{hempel}), together with the
Hurewicz theorem applied to the universal cover $\wwtilde{M}$.  This means
that $H^*(M)$ is isomorphic to the group cohomology $H^*(G)$, so algebraic
techniques can be applied.  In particular, construction of a chain
approximation to the diagonal (which we simply call a ``diagonal'') suffices
to determine the ring structure with arbitrary coefficients.

Most Seifert manifolds have infinite fundamental group: any Seifert
manifold with orbit surface not $S^2$ or $\RR P^2$, or having at least
four singular fibres, will have $G$ infinite.  Nevertheless, the relatively
small class of Seifert manifolds having finite fundamental group is extremely
important, indeed all known $3$--manifolds with finite fundamental group
are Seifert, and pending recent work of Perelman \cite{perelman}, 
Kleiner--Lott~\cite{kl}, Morgan--Tian \cite{mt} and
Cao--Zhu \cite{cz}, it seems very likely
there are no others.  These Seifert manifolds all arise from free orthogonal
actions of $G$ on $S^3$, and the resulting manifolds $M = S^3/G$, known as
spherical space forms, have been of great interest
to differential geometers since the nineteenth century;
see Clifford \cite{clifford}, Killing \cite{killing}, Klein \cite{klein} 
and the book of Wolf \cite{wolf}.  In this paper we attempt, in a certain sense,
to complete the aforementioned programme of Zieschang and his colleagues 
to the orientable Seifert manifolds with finite fundamental group, ie to the
spherical space forms.  (The nonorientable case has little interest here,
since a theorem of D\,B\,A\,Epstein \cite{epstein} asserts that $\ZZ_2$ is the only 
finite group that can be the fundamental group of a nonorientable $3$--manifold.)

It is important to note that, in contrast to the case where $G$ is infinite,
$M$ is no longer aspherical.  Thus, $H^*(M)$ and $H^*(G)$ are no longer
isomorphic; indeed by a classical theorem (see Cartan--Eilenberg \cite{ce}), 
$H^*(G)$ is now $4$--periodic.
The collection of all finite groups acting freely and orthogonally on $S^3$
is clearly listed by Milnor \cite{milnor}, based on earlier work of
Hopf \cite{hopf2} and Seifert--Threlfall \cite{st1}.  Ideally, for each
such group, one would like to have a $4$--periodic resolution ${\mathcal C}$
together with a contracting homotopy $s$ and a diagonal $\Delta$.  For
example, for the cyclic group $C_n$, this is done (here ${\mathcal C}$ is
$2$--periodic) in \cite{ce}.

In \fullref{sec:2}, we give some preliminaries about the groups involved and
about the cohomology of groups, also setting up necessary definitions and
notation.  The generalized quaternion groups $Q_{4n}$ are considered in
\fullref{sec:3}.  In this case, a $4$--periodic resolution was given in \cite{ce},
together with the somewhat cryptic statement ``the verification that the
homology groups are trivial involves some computations which will be omitted.''  
This verification was partially
done by Wall \cite{wall}, and is completely done here, ie, we give a
contracting homotopy $s$ for all $n \geq 1$.  A diagonal for $Q_{4n}$ was
first constructed by Shastri--Zvengrowski \cite{sz}.  The binary tetrahedral,
octahedral, and icosahedral groups (resp. $P_{24}$, $P_{48}$, $P_{120}$)
are discussed in \fullref{sec:4}.  Again, explicit $4$--periodic resolutions,
diagonals, and (for $P_{24}$, $P_{48}$) contracting homotopies are given.
The remaining two families of groups   $P'_{8\cdot3^k}$  and    $B_{2^k(2n+1)}$ 
are considered in \fullref{sec:5}.  Some concluding remarks, further questions,
and a brief discussion of applications, are given in \fullref{sec:6}.

For the most part, the results in this paper are given without proof.
This is partly because, once explicit formulae are found, 
the proofs are in general fairly routine
computations, but also because the verifications can often be quite
lengthy, eg the verification for the contracting homotopy and
diagonal map in \fullref{P:48} takes about 100 pages.  For full details,
see Tomoda \cite{tomoda}.

\section{Preliminaries} \label{sec:2}

In this section, we first discuss the groups that will be considered in 
the subsequent sections, namely the known finite fundamental groups of
$3$--manifolds.  In fact, every such group $G$ arises from a free orthogonal
action on $S^3$, with the resulting manifold $S^3/G$ an oriented Seifert 
manifold.  These groups were found in 1926 by Hopf \cite{hopf1}, and in
1931--33 by Seifert--Threlfall \cite{st1,st2}.  Further work in
1947 by Vincent \cite{vincent} considered the general case of free
orthogonal actions on any sphere (only the odd dimensional spheres are of
interest, since only $\ZZ_2$ can act freely on an even
dimensional sphere; cf Brown \cite{brown}).

The groups acting on $S^3$ were clearly listed (perhaps for the first
time) by Milnor in 1958 \cite{milnor}, as mentioned in
\fullref{sec:1}.  We denote them $C_n$, $Q_{4n}$, $n \geq 1$,
$P_{24}$, $P_{48}$, $P_{120}$, $B_{2^k(2n+1)}$, $k \geq 2$ and $n \geq
1$, $P'_{8\cdot3^k}$, $k \geq 1$, following Milnor's notation (except
that he denotes $B_{2^k(2n+1)}$ by $D_{2^k(2n+1)}$).  The direct
product of any of these groups with a cyclic group of relatively prime
order also acts freely and orthogonally on $S^3$.  In all cases, the
subscript denotes the order of the group, written $|G|$.  In Orlik's
1972 book \cite{orlik}, a considerably simplified derivation of this
list is given, but the shortest proof seems to be in a paper of
Hattori \cite{hattori} (in Japanese).  In the subsequent sections,
more details about each of these groups will be given, such as a
finite presentation and semidirect product structure.  From the work
of Milnor, Lee \cite{Lee} and Madsen--Thomas--Wall \cite{mtw}, there
remains the question concerning one other family of groups, $Q(8n, k,
l)$ (see \fullref{sec:6}), that could act freely on $S^3$ (or a
homotopy $S^3$).  Current work of Perelman \cite{perelman},
Kleiner--Lott \cite{kl}, Morgan--Tian \cite{mt} and Cao--Zhu \cite{cz}
will resolve this question (in the negative), as well as settle the
Poincar\'e conjecture and the geometrization conjecture for
$3$--manifolds.

We now briefly outline some of standard material about the cohomology of 
groups, following (chiefly) the book of Brown \cite{brown} as well as
other standard texts such as Adem--Milgram
\cite{am}, Benson \cite{benson1,benson2} and Cartan--Eilenberg \cite{ce}.  
Let $G$ be a finite group and $R = \ZZ G$ denote its integral
group ring.  An exact sequence ${\mathcal C}$ of projective (left) 
$R$--modules $C_j$, $j \geq 0$, and $R$--homomorphisms $d_j$, $j\geq 1$,
$${\mathcal C}: \hspace{.5cm}\cdots
\stackrel{d_{n+1}}{\lto} C_n\stackrel{d_n}{\lto}C_{n-1}\stackrel{d_{n-1}}{\lto} 
\cdots
\stackrel{d_2}{\lto} C_1\stackrel{d_1}{\lto}C_0
\stackrel{\ve}{\onto}\ZZ\to0 \hspace{.2cm},
$$
is called a projective resolution (in the subsequent sections, all resolutions will in
fact be free).  Here, $\ZZ$ has the trivial
$R$--module structure, and $\ve$ is called the augmentation.  It is also an
$R$--homomorphism, ie, $\ve(g\cdot x) = \ve(x)$, for all $g \in G$,
$x\in C_0$.  If $A$ is any (left) $R$--module, the cohomology of $G$ with
coefficients in $A$ is simply the cohomology of the cochain complex
$\hom_R({\mathcal C}, A)$, ie, $H^*(G; A) := H^*(\hom_R({\mathcal C}, A))$.

A contracting homotopy $s$ for ${\mathcal C}$ is a sequence of abelian group
homomorphisms $s_{-1}\co\ZZ\to C_0$ and $s_j\co C_j \to C_{j+1}$, $j \geq 0$ with
$\ve s_{-1} = 1_{\ZZ}$, $d_1s_0 + s_0\ve = 1_{C_0}$, 
$d_{j+1}s_j + s_{j-1}d_j = 1_{C_j}$, $j \geq 1$.  In general, $s_j$ is not
an $R$--homomorphism.  A contracting homotopy exists for any projective
resolution ${\mathcal C}$.

The chain complex $\ctimesc$ becomes a left
$R$--module via the diagonal action $g\cdot(x\otimes y) = gx\otimes gy$ for
$g \in G$, $x \in C_i$, $y \in C_j$, which is then extended by linearity 
to all of $R$.  A diagonal (strictly speaking, chain approximation to the 
diagonal) is an $R$--chain map
$\Delta\co {\mathcal C}\to\ctimesc$ such that

\begin{center}
\begin{picture}(80,40)
\put( 0,30){\makebox(0,0)[t]{\smash{$C_0$}}}
\put(60,30){\makebox(0,0)[t]{\smash{$C_0\otimes C_0$}}}
\put( 0, 0){\makebox(0,0)[t]{\smash{$\ZZ$}}}
\put(50, 0){\makebox(0,0)[t]{\smash{$\ZZ\iso\ZZ\otimes\ZZ$}}}

\put( 0,25){\vector(0,-1){15}}
\put(60,25){\vector(0,-1){15}}
\put( 6,32){\vector(1,0){33}}
\put( 5, 3){\vector(1,0){20}}

\put(-3,17){\makebox(0,0)[t]{\smash{$\ve$}}}
\put(72,17){\makebox(0,0)[t]{\smash{$\ve\otimes\ve$}}}
\put(25,35){\makebox(0,0)[t]{\smash{$\Delta_0$}}}
\put(14, 7){\makebox(0,0)[t]{\smash{$1_{\ZZ}$}}}

\end{picture}
\end{center}
commutes.

Using the resolution ${\mathcal C}$ and the diagonal map $\Delta$, the 
calculation of the cohomology $H^*(G; A)$ with coefficients in any 
$R$--module $A$ is quite routine, as well as the cup products when $A$ 
is an $R$--algebra.  In this paper, we content ourselves with a single illustration
of this process, in the proof of \fullref{ringP:482}, for $G = P_{48}$
with coefficients $\ZZ_2$.  
The calculation in all other cases can easily be reconstructed in the same
manner.

Exactness of a resolution can be proved by constructing a contracting
homotopy.  For a single finite group $G$, exactness can also be proved by
forgetting the $R$--module structure and simply showing exactness as a sequence 
of abelian groups, which is readily done with a computer 
(see Rotman \cite[p\,156]{rotman}).
The diagonal $\Delta$ can be used to determine the ring structure in
$H^*(G; A)$, where $A$ is any $R$--algebra.  Although $s$ and $\Delta$
always exist, finding either one explicitly can be a very demanding
calculation.  Once found, checking their required properties is relatively
routine, although often lengthy.

For a free resolution ${\mathcal C}$, a contracting homotopy $s$ can also
be used to produce a diagonal $\Delta$.  For example, following Handel \cite{handel},
we first define a contracting homotopy ${\tilde s}$ for $\ctimesc$ by
\begin{eqnarray*}
\tilde{s}_{-1} & = & s_{-1} \otimes s_{-1} \\
\tilde{s}_n\left(\sum_{i=0}^n(u_i\otimes v_{n-i})\right) & = &
\sum_{i=1}^ns_iu_i\otimes v_{n-i} + s_{-1}\ve(u_0)\otimes s_n(v_n),\ n \geq 0 \ ,
\end{eqnarray*}
where  $u_i\in C_i$, $v_{n-i}\in C_{n-i}$.

Then one defines $\Delta_n\co{\mathcal C} \to \ctimesc$ recursively on each free
generator $\rho_j$ of $C_n$ by
\begin{eqnarray*}
\Delta_0 & = & s_{-1}\ve\otimes s_{-1}\ve, \\
\Delta_n(\rho_j) & = & \tilde{s}_{n-1}\Delta_{n-1}d_n(\rho_j),
\end{eqnarray*}
and extends to all of $C_n$ by $R$--linearity.

\begin{Def}
A finite group $G$ is said to have \textit{periodic cohomology} of period $m$, 
if there exists a positive integer $m$ and a $u \in H^m(G; \ZZ) \iso \ZZ_{|G|}$ 
such that taking cup product with $u$ gives an isomorphism
$$
u\cup\rule{.3cm}{.4mm}\,\co H^l(G; A)\to H^{l+m}(G; A)
$$
for all $l \geq 1$ and for all $R$--modules $A$.

The element $u$ is called the periodicity class and $u\cup\rule{.3cm}{.4mm}$ is
called the periodicity isomorphism.  This definition can be given in
more elegant form, with the restriction $l \geq 1$ removed, using
Tate cohomology (see \cite[p\,260]{ce} and \cite[p\,153]{brown}).
\end{Def}

Any finite group $G$ acting freely on a sphere $S^{2n-1}$ will have 
$2n$--periodic cohomology, indeed, it will have a $2n$--periodic
resolution \cite{brown}.  Hence, the groups we study all have
$4$--periodic cohomology (with the cyclic groups $C_n$ being $2$--periodic).
The resolutions can be found by algebraic or geometric considerations.
Algebraically, it is advantageous to start with a balanced presentation
(same number of generators and relations) for $G$, then techniques of Fox
calculus will give $C_0$, $C_1$, and $C_2$ routinely.  For more
details, see \cite[Sections 2.3--2.4]{tomoda}.

The following sections consider the groups $Q_{4n}$, $P_{24}$, $P_{48}$,
$P_{120}$, $B_{2^k(2n+1)}$, $P'_{8\cdot 3^k}$.  Based mainly on the
dissertation of Tomoda \cite{tomoda}, we construct (as far as possible) a
$4$--periodic resolution ${\mathcal C}$ for each of these groups
together with a contracting homotopy
$s$ and a diagonal $\Delta$, as well as the cohomology ring $H^*(G; A)$ for
$A = \ZZ$, or  $A = \ZZ_p$ for a suitably selected prime $p$ (both as trivial $G$--modules).  
The cyclic groups $C_n$ are omitted since all this is completely done for
$C_n$ in \cite{ce}, and the corresponding orbit spaces $S^3/C_n$ are
the well known lens spaces ($\RR P^3$ for $n = 2$).  
We also omit the products $G\times C_n$
of any of the groups $G$ above with a cyclic group of relatively prime order,
since, for any groups $G_1$, $G_2$, 
$K(G_1\times G_2, 1) = K(G_1, 1) \times K(G_2, 1)$ implies that the cohomology
ring $H^*(G_1\times G_2)$ of the direct product of two groups can easily be
determined using the K\"unneth theorem.  Finally, for the associated
spherical space form $M = S^3/G$, note that $\pi_1(M) \iso G$ and
$\pi_j(M) \iso \pi_j(S^3)$, $j \geq 2$, from the homotopy exact sequence.
In particular, $\pi_2(M) = 0$, so by attaching cells to $M$ in dimensions
$4$ and higher, we see that the inclusion $i\co M \inc K(G, 1)$ embeds
$M$ as the $3$--skeleton of $K(G, 1)$.
It follows that $i^*\co H^l(G; A) \to H^l(M; A)$ is an isomorphism for $l \leq 2$
and, for $l = 3$, a monomorphism $H^3(G; A) \inj H^3(M; A) \iso A$.  Of course,
$H^l(M; A) = 0$ for $l \geq 4$.  Thus, it is not difficult to determine
$H^*(M; A)$ once $H^*(G; A)$ is known.
The following theorem briefly summarizes the results on the ring structures 
$H^*(M; A)$ for the spherical space forms $M = S^3/G$ (omitting the case $G$ 
cyclic, as mentioned above), with suitably chosen coefficient module(s) $A$.
The subscript of any cohomology class denotes its dimension.  
Since $H^l(M; A) = 0$ for $l > 3$, products of cohomology classes in total
dimension greater than $3$ are automatically $0$, so these relations are not
explicitly written in the polynomial rings below and are simply indicated by
the superscript ``$\star$.''
Further details, for each $G$, are given in the section devoted to that group.

\begin{Thm}
Using the notational conventions above, we have the following:
\begin{enumerate}
\item (cf~\fullref{genqua})\qua Let $M = S^3/Q_{4n}$, called a prism manifold 
\cite{orlik}. 

If $n$ is odd, then 
\begin{eqnarray*}
H^*(M; \ZZ_2) & \iso & \ZZ_2[\beta_1', \gamma_2', \delta_3]^{\star}
/(\left(\beta_1'\right)^2 = 0, \beta_1'\gamma_2' = \delta_3).
\end{eqnarray*}
If $n \equiv 0\ ({\rm mod\ } 4)$, then 
\begin{eqnarray*}
H^*(M; \ZZ_2) & \iso &\ZZ_2[\beta_1, \beta_1', \gamma_2, \gamma_2', \delta_3]^{\star}
\bigg/\left(\begin{array}{l}
  \beta_1^2 = \beta_1\beta_1' = \gamma_2, \left(\beta_1'\right)^2 = \gamma_2, \\
  \beta_1\gamma_2 = \beta_1\gamma_2' = \beta_1'\gamma_2' = \delta_3, \\
  \beta_1'\gamma_2 = 0
 \end{array}\right).
\end{eqnarray*}
If $n \equiv 2\ ({\rm mod\ } 4)$, then
\begin{eqnarray*}
H^*(M; \ZZ_2) & \iso & \ZZ_2[\beta_1, \beta_1', \gamma_2, \gamma_2', \delta_3]^{\star}
\bigg/\left(\begin{array}{l}
 \beta_1^2 = \gamma_2 + \gamma_2', \beta_1\beta_1' = \gamma_2', \\
 \left(\beta_1'\right)^2 = \gamma_2, \\
 \beta_1\gamma_2 = \beta_1\gamma_2' = \beta_1'\gamma_2' = \delta_3, \\
 \beta_1'\gamma_2 = 0
 \end{array}\right).
\end{eqnarray*}
\item (cf~\fullref{ringSMP:24})\qua Let $M = S^3/P_{24}$. 
\begin{eqnarray*}
H^*(M; {\mathbb Z}_3) & \iso & {\mathbb Z}_3[\beta_1, \gamma_2, \delta_3]^{\star}
/\left(\beta_1^2 = 0, \beta_1\gamma_2 = \delta_3\right).
\end{eqnarray*}
\item (cf~\fullref{ringSMP:48})\qua Let $M = S^3/P_{48}$. 
\begin{eqnarray*}
H^*(M; {\mathbb Z}_2) & \iso & {\mathbb Z}_2[\beta_1, \gamma_2, \delta_3]^{\star}
/\left(\beta_1^2 = \gamma_2, \beta_1\gamma_2 = \delta_3\right). \\
H^*(M; {\mathbb Z}_3) & \iso & {\mathbb Z}_3[\delta_3]^{\star}.
\end{eqnarray*}
\item (cf~\fullref{ringP:120})\qua Let $M = S^3/P_{120}$. 
The $3$--manifold $M$ is called the Poincar\'e homology sphere
and $H^l(M) = 0$ for all $l$ except $l = 0, 3$.  Thus, we have
$H^*(M; \ZZ) \iso \ZZ[\delta_3]^{\star}$ and $H^*(M; \ZZ_p) \iso \ZZ_p[\delta_3]^{\star}$.
\item (cf~\fullref{ringp'})\qua Let $M = S^3/P'_{8\cdot3^k}$. Then
\begin{align*}
H^*(M; {\mathbb Z}_3) & \iso  {\mathbb Z}_3[\beta_1, \gamma_2, \delta_3]^{\star}
/\left(\beta_1^2 = 0, \beta_1\gamma_2 = -\delta_3\right).
\\{\rm For} \ p \neq 3,  \ \ \ \ \ 
H^*(M; {\mathbb Z}_p)  &\iso  {\mathbb Z}_p[\delta_3]^{\star} . \end{align*}
\item (cf \fullref{thm5.6})\qua Let $M = S^3/B_{2^k(2n+1)}$, also called a 
prism manifold. Then
\begin{align*}
H^*(M; {\mathbb Z}_2) & \iso  \ZZ_2[\beta_1, \gamma_2, \delta_3 ]^{\star}
/\left(\beta_1^2 = 0, \beta_1\gamma_2 = \delta_3\right).
\\
 {\rm For} \ p \not= 2, \ \ \ \ \ 
H^*(M; {\mathbb Z}_p) &\iso {\mathbb Z}_p[\delta_3]^{\star}.\end{align*}
\end{enumerate}
\end{Thm}

\begin{Rem}
The above theorem includes all coefficients $\ZZ_p$, for those
primes $p$ of interest in each case (namely, $p$ divides the order
of $G_{ab} = H_1(M; \ZZ)$), as trivial $R$--modules.
For $\ZZ$ coefficients, see the corresponding section.
There are other possibilities for interesting (twisted) coefficients 
involving nontrivial $R$--modules; the authors hope to consider these 
in future work.
\end{Rem}

\section{Generalized quaternion groups} \label{sec:3}

In this section, we compute the ring structure of the cohomology of the 
generalized quaternion groups with $\ZZ$ and $\ZZ_2$ coefficients.
A presentation of the generalized quaternion groups is given by
$Q_{4n} = \<x, y\,|\, x^n = y^2, xyx=y\>$, for $n \geq 1$.  
One may also think of $Q_{4n}$ as a double cover of the dihedral group
$D_{2n} = \<\xi, \eta\,|\,\xi^n = \eta^2 = 1, \xi\eta\xi = \eta \>$,
using the exact sequence
$$
1 \to C_2 \stackrel{\vartriangleleft}{\inc} Q_{4n}
              \stackrel{p}{\onto}           D_{2n} \to 1
              \hspace{.2cm},
$$
where $C_2 = \{1, y^2\}$ is the centre of $Q_{4n}$ and $p(x) = \xi$,
$p(y) = \eta$.  This is related to the double cover 
$\Spin(3) \onto SO(3)$, indeed there is a commutative diagram
$$\begin{array}{ccccccccc}
1 & \to & C_2 & \stackrel{\vartriangleleft}{\inc} & Q_{4n}
              & \stackrel{p}{\onto}               & D_{2n} & \to & 1 \\
 & & \parallel & & \downarrow\subset & & \downarrow\subset & & \\
1 & \to & C_2 & \stackrel{\vartriangleleft}{\inc} & \Spin(3)
              & \stackrel{p}{\onto}               & SO(3)  & \to & 1 
              \hspace{.2cm}.
\end{array}$$
It is easy to show that 
$$
\left(Q_{4n}\right)_{ab} \iso \left\{\begin{array}{ll}
\ZZ_4,            & \mbox{ if $n$ is odd ,}     \\
\ZZ_2\oplus\ZZ_2, & \mbox{ if $n$ is even .}
\end{array}\right.
$$
A $4$--periodic resolution of $\ZZ$ over $R = \ZZ Q_{4n}$, $n \geq 1$,
will now be constructed (following Cartan--Eilenberg \cite{ce}).  First
define elements of $R$ as follows:
$$\begin{array}{lcll}
p_i & := & \sum_{k=0}^{i-1}x^k \hspace{.2cm}, & 0 \leq i \leq n \mbox{ with } p_0 := 0 \\
q_j & := & \sum_{k=0}^{j-1}y^k \hspace{.2cm}, & 0 \leq j \leq 4 \mbox{ with } q_0 := 0 \\
L   & := & p_n \hspace{.2cm},                              \\
N   & := & \sum_{g\in Q_{4n}}g \hspace{.2cm}.
\end{array}$$

\begin{Rem}
For any finite group $G$, following standard usage, the norm is written
$N := \sum_{g\in G}g\  \in \ZZ G$ (just as we did above for $Q_{4n}$).
\end{Rem}

\begin{Pro}
A resolution $\mathcal{C}$ for $Q_{4n}$ is given by:
$$\begin{array}{rlrl}
C_0 & =  \< a \>           & \ve(a) &= 1  \hspace{.2cm},\\
C_1 & =  \< b, b' \>   & d_1(b) &= (x-1)a \hspace{.2cm},\\
        &                                      & d_1(b') &= (y-1)a \hspace{.2cm},\\
C_2 & =  \< c, c' \>   & d_2(c) &= Lb - (y+1)b' \hspace{.2cm},\\
        &                                      & d_2(c') &= (xy+1)b + (x-1)b' \hspace{.2cm},\\
C_3 & =  \< d \>           & d_3(d) &= (x-1)c + (1-xy)c' \hspace{.2cm},\\
C_4 & =  \< a_4 \>         & d_4(a_4) &= Nd \hspace{.2cm}.
\end{array}$$
For any $n \geq 4$, we define $C_n \approx C_{n-4}$ with appropriate subscripts, similarly
$d_n$ is defined in the obvious way from $d_{n-4}$ (note that in the above resolution,
strictly speaking, $a=a_0$,  $b=b_1$, etc).
\end{Pro}

The resolution above is given in \cite{ce} without proof.
Wall showed in \cite{wall} that the chain complex ${\mathcal C}$ above
is a resolution for $n$ even via representation theory.
The following contracting homotopy verifies directly that the chain complex
${\mathcal C}$ above is indeed a resolution of $\ZZ$ over $\ZZ Q_{4n}$,
for all $n \geq 1$, thus completing the claim of Cartan--Eilenberg and the
work of Wall.

\begin{Pro}
Let $0 \leq i \leq n-1$ and $0 \leq j \leq 3$.  Then a contracting homotopy $s$ 
for ${\mathcal C}$ is given by: 
$$\begin{array}{rlrl}
s_{-1}(1)                & =  a \hspace{.2cm}, \\
 \\
s_0(x^iy^ja)         & =  p_ib + x^iq_jb' \hspace{.2cm}, \\
 \\
s_1(x^ib)                & =  0 \hspace{.2cm},                                     && 0\leq i\leq n-2 \\
s_1(x^{n-1}b)        & =  c \hspace{.2cm},                                      \\
s_1(yb)                  & =  (x^{n-1}-1+x^{n-1}y)c + (y -x^{n-1}yL)c' \hspace{.2cm}, 
                                                                                         \\
s_1(x^iyb)           & =  x^{i-1}c' \hspace{.2cm},                          && 1\leq i\leq n-1 \\
s_1(x^iy^2b)         & =  x^i(x-1)c \hspace{.2cm},                          && 0\leq i\leq n-2 \\
 \end{array}$$$$\begin{array}{rlrl}
s_1(x^{n-1}y^2b) & =  -yc + (yL+y^3-x^{n-1})c' \hspace{.2cm},  \\
s_1(y^3b)                & =  -c + x^{n-1}c' \hspace{.2cm},                 \\
s_1(x^iy^3b)         & =  x^{i-1}(1-xy+y^2)c' \hspace{.2cm},    && 1\leq i\leq n-1 \\
 \\
s_1(x^iy^jb')        & =  0 \hspace{.2cm},                                          && 0\leq j\leq2 \\
s_1(x^iy^3b')        & =  -x^i(y+1)c + x^iyLc' \hspace{.2cm},  \\
  \\ 
s_2(x^ic)                   & =  0 \hspace{.2cm},                     \\
s_2(x^iyc)              & =  (p_{n-i-1}x^{i+1}y - p_{n-1})d \hspace{.2cm},  \\
s_2(x^iy^2c)            & =  (x^{n-1}p_{i+1} + p_ixy)d \hspace{.2cm},  \\
s_2(x^iy^3c)            & =  (x^{i+1}Ly^3 - p_i)d \hspace{.2cm},  \\
 \\
s_2(x^ic')              & =  0 \hspace{.2cm},                                       && 0\leq i\leq n-1 \\
s_2(yc')                    & =  0  \hspace{.2cm},                                      \\
s_2(x^iyc')             & =  -x^{i-1}d \hspace{.2cm},                   && 1\leq i\leq n-1 \\
s_2(x^iy^2c')           & =  -x^id \hspace{.2cm},                           && 0\leq i\leq n-2 \\
s_2(x^{n-1}y^2c')   & =  (x^{n-1} + Ly^2 + p_{n-1}xy)d \hspace{.2cm},  \\
s_2(y^3c')              & =  0  \hspace{.2cm},                                       \\
s_2(x^iy^3c')           & =  x^{i-1}(xy-1)d \hspace{.2cm},                 && 1\leq i\leq n-1 
 \\ \\ 
s_3(y^3d)                   & =  a_4  \hspace{.2cm},                   \\
s_3(x^iy^3d)            & =  0 \hspace{.2cm},                      && \hbox{ otherwise.}
\end{array}$$
The remaining $s_n$, for $n \geq 4$, are then defined by periodic
extension, for example, $s_4(x^iy^ja_4) = p_ib_5 + x^iq_jb'_5$, etc.
\end{Pro}

As mentioned in the Introduction, the proofs for this proposition and
most of the following ones are not given here, for full details, see Tomoda
\cite{tomoda}.
The following defines a diagonal map
$\Delta\co{\mathcal C} \to {\mathcal C} \otimes {\mathcal C}$
for ${\mathcal C}$ through dimension $4$.  We remark that the contracting
homotopy $s$ extends to higher dimensions by periodicity, as noted above,
but this is not true for $\Delta$.

\begin{Pro}
A diagonal map $\Delta$ for ${\mathcal C}$ is given by:
$$\begin{array}{lcl}
\Delta_0(a)  & = & a\otimes a \hspace{.2cm}, \\
\Delta_1(b)  & = & b\otimes xa + a\otimes b \hspace{.2cm}, \\
\Delta_1(b') & = & b'\otimes ya + a\otimes b' \hspace{.2cm}, \\
\Delta_2(c)  & = &
  c\otimes y^2a + \sum_{i=0}^{n-1}(p_ib\otimes x^ib) 
  + a\otimes c - b'\otimes yb' \hspace{.2cm}, \\
\Delta_2(c') & = &
  c'\otimes ya + b\otimes xyb + xb'\otimes xyb + a\otimes c'
    + b\otimes xb' \hspace{.2cm}, \\
     \end{array}$$$$\begin{array}{lcl}
\Delta_3(d)  & = &
    c\otimes y^2b + b\otimes xc + d\otimes xy^2a - c'\otimes y^2b \\
             &   &  
  - b\otimes xyc' - xb'\otimes xyc' + a\otimes d - c'\otimes yb'  \hspace{.2cm},  \\
\Delta(a_4) 
& = & 
  a\otimes a_4 
  +\sum_{i=0}^{n-1}\sum_{j=0}^3(p_ib\otimes x^iy^jd) 
\\
  &&+ \sum_{i=0}^{n-1}\sum_{j=0}^3(x^iq_jb'\otimes x^iy^jd)   \\
 &   &
 + c\otimes y^2c + (x^{n-1}-1+x^{n-1}y)c\otimes x^{n-1}y^3c \\
 &   & 
 + \sum_{i=0}^{n-2}(x^i(x-1)c\otimes x^{i+1}y^2c)
 - yc\otimes c - c\otimes x^{n-1}yc                         \\
 &   & 
 - c\otimes y^3c' - (x^{n-1}-1+x^{n-1}y)c\otimes x^{-1}y^2c' \\
 &   & 
 - \sum_{i=0}^{n-2}(x^i(x-1)c\otimes x^{i+1}y^3c') + yc\otimes yc' \\
 &   & 
 + c\otimes x^{n-1}y^2c' + \sum_{i=0}^{n-1}(x^i(y+1)c\otimes x^iyc') \\
 &   & 
 + (y -x^{n-1}yL)c'\otimes x^{-1}yc
 + \sum_{i=1}^{n-1}(x^{i-1}c'\otimes x^{i-1}yc)             \\
 &   & 
 + (yL+y^3-x^{n-1})c'\otimes x^ny^2c + x^{n-1}c'\otimes x^{-1}y^3c \\
 &   & 
 + \sum_{i=1}^{n-1}(x^{i-1}(1-xy+y^2)c'\otimes x^{i-1}y^3c) \\
 &   & 
 - (y -x^{n-1}yL)c'\otimes x^{-1}y^2c'
 - \sum_{i=1}^{n-1}(x^{i-1}c'\otimes x^{i-1}y^2c')          \\
 &   & 
 - (yL+y^3-x^{n-1})c'\otimes x^ny^3c' - x^{n-1}c'\otimes x^{-1}c' \\
 &   & 
 - \sum_{i=1}^{n-1}(x^{i-1}(1-xy+y^2)c'\otimes x^{i-1}c')
 - \sum_{i=0}^{n-1}(x^iyLc'\otimes x^iyc')                  \\
 &   & 
 + \sum_{i=0}^{n-1}((p_{n-i-1}x^{i+1}y - p_{n-1})d\otimes x^iy^3b) \\
 &   &
 + \sum_{i=0}^{n-1}((x^{n-1}p_{i+1} + p_ixy)d\otimes x^ib) \\
 &   & 
 + \sum_{i=0}^{n-1}((x^{i+1}Ly^3 - p_i)d\otimes x^iyb)      \\
 &   & 
 - \sum_{i=1}^{n-1}(-x^{i-1}d\otimes x^iy^3b) 
 - \sum_{i=0}^{n-2}(-x^id\otimes x^ib)                      \\
 &   & 
 - (x^{n-1} + Ly^2 + p_{n-1}xy)d\otimes x^{n-1}b 
\\&& - \sum_{i=1}^{n-1}(x^{i-1}(xy-1)d\otimes x^iyb)            
\\
\hspace{1.3cm} &   & 
 - \sum_{i=1}^{n-1}(-x^{i-1}d\otimes x^iy^2b')
 - \sum_{i=0}^{n-2}(-x^id\otimes x^iy^3b')                  \\
 &   & 
 - (x^{n-1} + Ly^2 + p_{n-1}xy)d\otimes x^iy^3b' 
\\ &&- \sum_{i=1}^{n-1}(x^{i-1}(xy-1)d\otimes x^ib')            \\
 &   & 
 + a_4\otimes x^{n-1}y^3a  \hspace{.2cm}. 
\end{array}$$
\end{Pro}

\begin{Pro}
The cohomology groups of the generalized quaternion group $Q_{4n}$, 
for $n \geq 1$, are given by:
\begin{eqnarray*}
H^l(Q_{4n}; {\mathbb Z})      & = &
\left\{\begin{array}{ll}
  {\mathbb Z} ,   &\mbox{ if $l = 0$, } \\
  0 ,             &\mbox{ if $l \equiv 1$ {\rm mod }$4$,} \\
  \left\{\begin{array}{l}
    {\mathbb Z}_2 \oplus {\mathbb Z}_2 , \\
    {\mathbb Z}_4 , 
  \end{array}\right. & \begin{array}{l}
       \mbox{if $l \equiv 2$ {\rm mod }$4$ and $n$ even,} \\
    \mbox{if $l \equiv 2$ {\rm mod }$4$ and $n$ odd,} 
  \end{array}
  \\
  0 ,              &\mbox{ if $l \equiv 3$ {\rm mod }$4$,} \\
  {\mathbb Z}_{4n} ,  
                            &\mbox{ if $l \equiv 0$ {\rm mod }$4$ and $l > 0$.} \\
\end{array}\right.
\end{eqnarray*}
\end{Pro}

\begin{Thm} \label{ringgenquat}
The cohomology ring
$H^*(Q_{4n}; {\mathbb Z})$ has the following presentation:
$$
      H^*(Q_{4n}; {\mathbb Z}) 
\iso 
\left\{\begin{array}{ll}
  {\mathbb Z}[\gamma_2, \gamma_2', \alpha_4]\big/
   \left(\begin{array}{l}
   2\gamma_2 = 2\gamma_2' = 0 = 4n\alpha_4,\\
   \gamma_2^2 = 0, \gamma_2\gamma_2' = \gamma_2'^2 = 2n\alpha_4 
   \end{array}\right) ,  
  & \mbox{if } n = 4m, \\
  {\mathbb Z}[\gamma_2', \alpha_4]/
   \left(4\gamma_2' = 0 = 4n\alpha_4, \gamma_2'^2 = n\alpha_4 \right) ,  
  & \mbox{if } n = 4m + 1, \\
  {\mathbb Z}[\gamma_2, \gamma_2', \alpha_4]\big/
   \left(\begin{array}{l} 
   2\gamma_2 = 2\gamma_2' = 0 = 4n\alpha_4, \\
   \gamma_2^2 = 0 = \gamma_2'^2, \gamma_2\gamma_2' = 2n\alpha_4 
   \end{array}\right) , 
  & \mbox{if } n = 4m + 2, \\
  {\mathbb Z}[\gamma_2', \alpha_4]/
    \left( 4\gamma_2' = 0 = 4n\alpha_4, \gamma_2'^2 = 3n\alpha_4 \right) ,  
  & \mbox{if } n = 4m + 3 \, .
\end{array}\right.
$$
\end{Thm}

\begin{Pro}
The cohomology groups of the generalized quaternion group $Q_{4n}$
with $\ZZ_2$ coefficients, for $n \geq 1$, are given by:
\begin{eqnarray*}
H^l(Q_{4n}; \ZZ_2)      & = &
\left\{\begin{array}{l}
  {\mathbb Z}_2 \hspace{.2cm},  \hspace{1.7cm}
      \mbox{ if $l \equiv 0,\ 1$ {\rm mod }$4$, }             \\
  \left\{\begin{array}{ll}
    {\mathbb Z}_2 \oplus {\mathbb Z}_2 \hspace{.2cm}, 
    & \mbox{ if $l \equiv 2,\ 3$ {\rm mod }$4$ and $n$ even,} \\
    {\mathbb Z}_2 \hspace{.2cm}, 
    & \mbox{ if $l \equiv 2,\ 3$ {\rm mod }$4$ and $n$ odd.}  \\
  \end{array}\right. \\
\end{array}\right.
\end{eqnarray*}
\end{Pro}

\begin{Thm}
For $n \equiv 0\ ({\rm mod}\ 4)$, the cohomology ring $H^*(Q_{4n}; \ZZ_2)$ 
is given by:
$$ 
H^*(Q_{4n}; \ZZ_2) \iso
  {\mathbb Z}_2[\beta_1, \beta_1', \gamma_2, \gamma_2', \delta_3, \alpha_4]\bigg/
    \left( \begin{array}{l}
           \beta_1^2 = \gamma_2' = \beta_1\beta_1',
           \left(\beta_1'\right)^2 = \gamma_2, \\
           \beta_1\gamma_2 = \beta_1\gamma_2' = \beta_1'\gamma_2' = \delta_3, \\
           \beta_1'\gamma_2 = 0 , \\
           \gamma_2^2 = \left(\gamma_2'\right)^2 = \gamma_2\gamma_2' = 0
           \end{array}
    \right), 
$$
and for $n \equiv 2\ ({\rm mod}\ 4)$,
$$
H^*(Q_{4n}; \ZZ_2) \iso 
  {\mathbb Z}_2[\beta_1, \beta_1', \gamma_2, \gamma_2', \delta_3, \alpha_4]\bigg/
    \left( \begin{array}{l}
           \beta_1^2 = \gamma_2+\gamma_2', \beta_1\beta_1' = \gamma_2', \\
           \left(\beta_1'\right)^2 = \gamma_2, \\
           \beta_1\gamma_2 = \beta_1\gamma_2' = \beta_1'\gamma_2' = \delta_3, \\
           \beta_1'\gamma_2 = 0 , \\
           \gamma_2^2 = \left(\gamma_2'\right)^2 = \gamma_2\gamma_2' = 0
           \end{array}
    \right). 
$$
For $n$ odd, the cohomology ring $H^*(Q_{4n}; \ZZ_2)$ is given by:
$$
H^*(Q_{4n}; \ZZ_2) \iso
  {\mathbb Z}_2[\beta_1', \gamma_2', \delta_3, \alpha_4]\bigg/
    \left( \begin{array}{l}
           \left(\beta_1'\right)^2 = 0, 
           \beta_1'\gamma_2' = \delta_3, \\ 
           \beta_1'\delta_3 = 0, 
           \left(\gamma_2'\right)^2 = \alpha_4 \\
           \end{array}
    \right). 
$$
\end{Thm}

\begin{Cor}\label{genqua}
Let $M = S^3/Q_{4n}$.  Then the following holds:
\begin{enumerate}
\item

$
     H^*(M; {\mathbb Z}) 
\iso 
\left\{\begin{array}{ll}
  {\mathbb Z}[\gamma_2, \gamma_2', \delta_3]^{\star}/
  (2\gamma_2 = 2\gamma_2' = 0) \hspace{.2cm},  
  & \mbox{ if\ $n$\ is\ even }, \\
  {\mathbb Z}[\gamma_2', \delta_3]^{\star}/
  (4\gamma_2' = 0) \hspace{.2cm},  
  & \mbox{ if } n \mbox{\ is\ odd}.
\end{array}\right.
$
\item When $n \equiv 0 ({\rm mod}\ 4)$, \\
$H^*(M; \ZZ_2) \iso \ZZ[\beta_1, \beta_1', \gamma_2, \gamma_2', \delta_3]^{\star}
\bigg/ \left(\begin{array}{l}
       \beta_1^2 = \gamma_2\ = \beta_1\beta_1',
       \left(\beta_1'\right)^2 = \gamma_2, \\
       \beta_1\gamma_2 = \beta_1\gamma_2' = \beta_1'\gamma_2' = \delta_3, \\
       \beta_1'\gamma_2 = 0 \\
 \end{array} \right)\, .$
\item When $n \equiv 2 ({\rm mod}\ 4)$, 
$$H^*(M; \ZZ_2) \iso \ZZ[\beta_1, \beta_1', \gamma_2, \gamma_2', \delta_3]^{\star}
\bigg/ \left(\begin{array}{l}
       \beta_1^2 = \gamma_2 + \gamma_2', \beta_1\beta_1' = \gamma_2', \\
       \left(\beta_1'\right)^2 = \gamma_2, \beta_1'\gamma_2 = 0, \\
       \beta_1\gamma_2 = \beta_1\gamma_2' = \beta_1'\gamma_2' = \delta_3
 \end{array} \right) \, .$$
\item When $n$ is odd, 
$$H^*(M; \ZZ_2) \iso \ZZ[\beta_1', \gamma_2', \delta_3]^{\star}
/ \left(\begin{array}{l}
       \left(\beta_1'\right)^2 = 0,
       \beta_1'\gamma_2' = \delta_3
 \end{array} \right) \, .$$
\end{enumerate}
\end{Cor}

\section{Binary groups} \label{sec:4}

In this section, we consider double covers (under the $2$--fold covering
$\Spin(3)\onto SO(3)$) of the tetrahedral, octahedral, and icosahedral
groups, called respectively the binary tetrahedral, binary octahedral,
and binary icosahedral groups.  The generalized quaternion groups
$Q_{4n}$, considered in \fullref{sec:3}, could also be thought of as ``binary
 dihedral groups.''

\subsection{Binary tetrahedral group} \label{P:24}
The binary tetrahedral group $P_{24}$ can be considered as a double
cover of the group of rotational symmetries ${\mathfrak A}_4$ of a
regular tetrahedron (${\mathfrak A}_4$ is the alternating group on the $4$
symbols $\{1, 2, 3, 4\}$).  Thus, there is a commutative diagram of short
exact sequences
$$\begin{array}{ccccccccc}
1 & \to & C_2 & \stackrel{\vartriangleleft}{\inc} & P_{24}
              & \stackrel{p}{\onto}               & {\mathfrak A}_4 & \to & 1 \\
 & & \parallel & & \downarrow\subset & & \downarrow\subset & & \\
1 & \to & C_2 & \stackrel{\vartriangleleft}{\inc} & \Spin(3)
              & \stackrel{p}{\onto}               & SO(3)           & \to & 1 
              \hspace{.2cm}.
\end{array}$$
Following the book of Coxeter--Moser \cite{cm}, we use the balanced presentation
$P_{24} = \< S, T \,|\, STS = T^2, TST = S^2 \>$.  It is easy to see that
$z := (ST)^2 = T^3 = (TS)^2 = S^3$, and this element is central.  Then
$C_2 = \{1, z\}$ is the centre of $P_{24}$.  The homomorphism $p$ is given 
by $p(S) = (1\,2\,3) \in {\mathfrak A}_4$, 
$p(T) = (1\,2\,4) \in {\mathfrak A}_4$
(note that $p$ is not unique) .
It is easy to show $\left(P_{24}\right)_{ab} \iso \ZZ_3$.

Other common presentations of $P_{24}$ are
$\< x, y \,|\, x^2 = (xy)^3 = y^3, x^4 = 1 \>$ and
$\< x, y \,|\, x^2 = (xy)^3 = y^{-3} \>$.  The equivalence can be established
using $x = ST$ and $y = T^{-1}$.

\begin{Pro}
A resolution $\mathcal{C}$ for $P_{24}$ is given by:
$$\begin{array}{rlrl}
C_0 & =  \<a\>         & \varepsilon(a) &= 1 \hspace{.2cm}, \\
C_1 & =  \<b, b'\> & d_1(b) &= (S-1)a \hspace{.2cm}, \\
        &                              & d_1(b') &= (T-1)a \hspace{.2cm}, \\
C_2 & =  \<c, c'\> & d_2(c) &= (T-S-1)b + (1+TS)b' \hspace{.2cm}, \\
        &                              & d_2(c') &= (1+ST)b + (S-T-1)b' \hspace{.2cm}, \\
C_3 & =  \<d\>         & d_3(d) &= (S-1)c + (T-1)c' \hspace{.2cm}, \\
C_4 & =  \<a_4\>   & d_4(a_4) &= Nd \hspace{.2cm}.
\end{array}$$
For any $n \geq 4$, we define $C_n \approx C_{n-4}$ with appropriate subscripts.
\end{Pro}

We now define a contracting homotopy for this resolution.

\begin{Pro}
A contracting homotopy $s$ for the resolution ${\mathcal C}$ over
$\ZZ P_{24}$ above is given by:
$$
\begin{array}{rlrl}
s_{-1}(1) & =  a \hspace{.2cm},&\\
\\
s_0(a) & = 0 \hspace{.2cm},& \qquad\ \ s_0(TSa) & =  Tb + b' \hspace{.2cm},  \\

s_0(Sa) & =  b \hspace{.2cm},  &\qquad\ \ s_0(S^2a) & =  (1 + S)b \hspace{.2cm},  \\
s_0(Ta) & =  b' \hspace{.2cm},  & \qquad\ \ s_0(T^2a) & =  (1 + T)b' \hspace{.2cm},  \\
s_0(STa) & =  Sb'+b \hspace{.2cm},  &\qquad\quad s_0(ST^2a) & =  b + S(1 + T)b' \hspace{.2cm}, 
\end{array}
$$
$$
\begin{array}{rl}
s_0(TS^2a) & =  T(1 + S)b + b' \hspace{.2cm},  \\
s_0(S^2Ta) & =  (1 + S)b + S^2b' \hspace{.2cm},  \\
s_0(T^2Sa) & =  T^2b + (1 + T)b' \hspace{.2cm},  \\
s_0(ST^2Sa) & =  (1+ ST^2)b + S(1 + T)b' \hspace{.2cm},  \\
\\
s_0(z a) & =  (1 + ST)b + (S + T^2)b' \hspace{.2cm},   \\
s_0(z Sa) & =  z b + (1 + ST)b + (S + T^2)b' \hspace{.2cm},   \\
s_0(z Ta) & =  z b' + (1 + ST)b + (S + T^2)b' \hspace{.2cm},   \\
s_0(z STa) & =  z(b + Sb') + (1 + ST)b + (S + T^2)b' \hspace{.2cm},   \\
s_0(z TSa) & =  z(Tb + b') + (1 + ST)b + (S + T^2)b' \hspace{.2cm},   \\
s_0(z S^2a) & =  z(S + 1)b + (1 + ST)b + (S + T^2)b' \hspace{.2cm},   \\
s_0(z T^2a) & =  z(T + 1)b' + (1 + ST)b + (S + T^2)b' \hspace{.2cm},   \\
s_0(z ST^2a) & =  z(b + S(1 + T)b') + (1 + ST)b + (S + T^2)b' \hspace{.2cm},   \\
s_0(z TS^2a) & =  z(T(1 + S)b + b') + (1 + ST)b + (S + T^2)b' \hspace{.2cm},   \\
\end{array}
$$
$$
\begin{array}{rl}
s_0(z S^2Ta) & =  z((1 + S)b + S^2b') + (1 + ST)b + (S + T^2)b' \hspace{.2cm},   \\
s_0(z T^2Sa) & =  z(T^2b + (1 + T)b') + (1 + ST)b + (S + T^2)b' \hspace{.2cm},   \\
s_0(z ST^2Sa) & =  z((1+ ST^2)b + S(1 + T)b') + (1 + ST)b + (S + T^2)b' \hspace{.2cm},   \\
 &  \\
s_1(b) & =  0 \hspace{.2cm},  \\
s_1(Sb) & =  0 \hspace{.2cm},  \\
s_1(Tb) & =  0 \hspace{.2cm},  \\
s_1(STb) & =  c' \hspace{.2cm},  \\
s_1(TSb) & =  0 \hspace{.2cm},  \\
s_1(S^2b) & =  -Sc \hspace{.2cm},  \\
s_1(T^2b) & =  0 \hspace{.2cm},  \\
s_1(ST^2b) & =  0 \hspace{.2cm},  \\
s_1(TS^2b) & =  -TSc + (T - 1)c' \hspace{.2cm},  \\
s_1(S^2Tb) & =  Sc' \hspace{.2cm},  \\
s_1(T^2Sb) & =  -T^2c - c' \hspace{.2cm},   \\
s_1(ST^2Sb) & =  - ST^2c - STc' \hspace{.2cm},  \\
 & \\
s_1(z b) & =  0 \hspace{.2cm},  \\
s_1(z Sb) & =  0 \hspace{.2cm},  \\
s_1(z Tb) & =  0 \hspace{.2cm},  \\
s_1(z STb) & =  z c' \hspace{.2cm},  \\
s_1(z TSb) & =  0 \hspace{.2cm},  \\
s_1(z S^2b) & =  (T + T^2 + z S^2T)c + (1 + T^2S + z S + z S^2)c' \hspace{.2cm},  
\\ 
s_1(z T^2b) & =  0 \hspace{.2cm},  \\
s_1(z ST^2b) & =  0 \hspace{.2cm},  \\
s_1(z TS^2b) & =  (S^2 + T^2 + S^2T)c + (1 + S^2 + ST^2 + z TS^2)c' \hspace{.2cm},  \\
s_1(z S^2Tb) & =  z Sc' \hspace{.2cm},  \\
s_1(z T^2Sb) & =  (S + TS^2 + z ST^2S)c + (S + S^2 + TS^2)c' \hspace{.2cm},  \\
s_1(z ST^2Sb) & =  (S + TS + T^2S)c + (TS + T^2S + z)c' \hspace{.2cm},  \\
\end{array}
$$
$$
\begin{array}{rlrl}
s_1(b') & =  0 \hspace{.2cm},  & s_1(TS^2b') & =  -(T + T^2)c - (1 + T^2S)c' \hspace{.2cm},  \\
s_1(Sb') & =  0 \hspace{.2cm},  & s_1(S^2Tb') & =  -Sc - S^2c' \hspace{.2cm},  \\
s_1(Tb') & =  0 \hspace{.2cm},  & s_1(T^2Sb') & =  Tc \hspace{.2cm},  \\
s_1(STb') & =  0 \hspace{.2cm},  & s_1(ST^2Sb') & =  STc + c' \hspace{.2cm},  \\
s_1(TSb') & =  c \hspace{.2cm},  \\
s_1(S^2b') & =  0 \hspace{.2cm},  \\
s_1(T^2b') & =  -c' \hspace{.2cm},  \\
s_1(ST^2b') & =  -STc' \hspace{.2cm},  \\
\end{array}
$$
$$
\begin{array}{rl}
\mskip96mu s_1(z b') & =  0 \hspace{.2cm},  \\
s_1(z Sb') & =  0 \hspace{.2cm},  \\
s_1(z Tb') & =  0 \hspace{.2cm},  \\
s_1(z STb') & =  0 \hspace{.2cm},  \\
s_1(z TSb') & =  z c \hspace{.2cm},  \\
s_1(z S^2b') & =  0 \hspace{.2cm},  \\
s_1(z T^2b') & =  (S + TS + T^2S + z S^2T)c
                   + (T^2S + z ST + z S^2T)c' \hspace{.2cm},  \\
s_1(z ST^2b') & =  (S + TS + T^2S + z S^2T)c
                    + (T^2S + z + z S^2T)c' \hspace{.2cm},  \\
s_1(z TS^2b') & =  (S^2T + z)c + (ST + S^2T)c' \hspace{.2cm}, \\
s_1(z S^2Tb') & =  (T + T^2 + z S^2T)c + (1 + T^2S + z S)c' \hspace{.2cm},   \\
s_1(z T^2Sb') & =  z Tc \hspace{.2cm},   \\
s_1(z ST^2Sb') & =  z (STc + c') \hspace{.2cm},  \\
\end{array}
$$
$$
\begin{array}{rlrl}
s_2(c) & =  0 \hspace{.2cm},  &\qquad\qquad\ \  s_2(TSc) & =  d \hspace{.2cm},  \\
s_2(Sc) & =  0 \hspace{.2cm},  & s_2(S^2c) & =  Sd \hspace{.2cm},  \\
s_2(Tc) & =  0 \hspace{.2cm},  & s_2(T^2c) & =  0 \hspace{.2cm},  \\
s_2(STc) & =  0 \hspace{.2cm},  & s_2(ST^2c) & =  0 \hspace{.2cm},  \\
s_2(TS^2c) & =  T(-1 + S - T)d \hspace{.2cm},  \\
s_2(S^2Tc) & =  -S(1 + S)d \hspace{.2cm},  \\
s_2(T^2Sc) & =  T(1 + T)d \hspace{.2cm},  \\
s_2(ST^2Sc) & =  ST(1 + T)d \hspace{.2cm},  \\
\end{array}
$$
$$
\begin{array}{rl}
s_2(z c) & =  0 \hspace{.2cm},  \\
s_2(z Sc) & =  -(1 + T + T^2 + z S + z S^2)d \hspace{.2cm},  \\
s_2(z Tc) & =  0 \hspace{.2cm},  \\
s_2(z STc) & =  0 \hspace{.2cm},  \\
s_2(z TSc) & =  -(S^2 + ST^2 + S^2T + ST^2S + z TS + z TS^2)d \hspace{.2cm},  \\
s_2(z S^2c) & =  -(1 + T + T^2 + z S^2)d \hspace{.2cm},  \\
s_2(z T^2c) & =  (TS + T^2S + z ST + z ST^2 + z S^2T)d \hspace{.2cm},  \\
s_2(z ST^2c) & =  z S^2Td \hspace{.2cm},  \\
s_2(z TS^2c) & =  -(S^2 + ST^2 + ST^2S + z TS^2)d \hspace{.2cm},  \\
s_2(z S^2Tc) & =  0 \hspace{.2cm},  \\
s_2(z T^2Sc) & =  -(S + ST + TS^2 + z T^2S + z ST^2S)d \hspace{.2cm},  \\
s_2(z ST^2Sc) & =  -(TS + T^2S)d \hspace{.2cm},  \\
\\
s_2(c') & =  0 \hspace{.2cm},  \\
s_2(Sc') & =  0 \hspace{.2cm},  \\
s_2(Tc') & =  d \hspace{.2cm},  \\
s_2(STc') & =  0 \hspace{.2cm},  \\
 \end{array}
$$
$$
\begin{array}{rl}
s_2(TSc') & =  -(1 + T + T^2)d \hspace{.2cm},  \\
s_2(S^2c') & =  0 \hspace{.2cm},  \\
s_2(T^2c') & =  Td \hspace{.2cm},  \\
s_2(ST^2c') & =  STd \hspace{.2cm},  \\
s_2(TS^2c') & =  T(1 + T + TS)d \hspace{.2cm},  \\
s_2(S^2Tc') & =  S(1 + S)d \hspace{.2cm},   \\
s_2(T^2Sc') & =  0 \hspace{.2cm},  \\
s_2(ST^2Sc') & =  -ST(1 + T + TS)d \hspace{.2cm},  \\
\\
s_2(z c') & =  0 \hspace{.2cm},  \\
s_2(z Sc') & =  0 \hspace{.2cm},  \\
s_2(z Tc') & =  (1 + T + T^2 + z(1 + S + S^2))d \hspace{.2cm},  \\
s_2(z STc') & =  0 \hspace{.2cm},  \\
s_2(z TSc') & =  S^2Td \hspace{.2cm},  \\
s_2(z S^2c') & =  0 \hspace{.2cm},  \\
s_2(z T^2c') & =  (1 + T + S^2 + T^2 + ST^2 + S^2T + ST^2S \\
                  & + z + z S + z T + z TS + z S^2 + z TS^2)d \hspace{.2cm},  \\
s_2(z ST^2c') & =  -(TS + T^2S + z ST^2 + z S^2T)d \hspace{.2cm},  \\
s_2(z TS^2c') & =  (-ST + S^2)d \hspace{.2cm},  \\
s_2(z S^2Tc') & =  -(1 + T + T^2)d \hspace{.2cm},  \\
s_2(z T^2Sc') & =  (S + S^2 + TS^2 + z ST^2S)d \hspace{.2cm},  \\
s_2(z ST^2Sc') & =  (TS + T^2S + TS^2)d \hspace{.2cm},  \\
&  \\
s_3(z T^2d) & =  a_4 \hspace{.2cm}. 
\end{array}$$
\end{Pro}

\begin{Pro}
For the given resolution of
${\mathcal C}\stackrel{\ve}{\onto}{\mathbb Z}$
over ${\mathbb Z}P_{24}$, a diagonal map
$\Delta\co{\mathcal C} \to {\mathcal C}\otimes{\mathcal C}$,
through dimension $4$, is given by:
$$\begin{array}{lcl}
\Delta_0(a)  & = &
a \otimes a \hspace{.2cm},  \\
\Delta_1(b)  & = &
b \otimes Sa + a \otimes b \hspace{.2cm},  \\
\Delta_1(b') & = &
b' \otimes Ta + a \otimes b'  \hspace{.2cm}, \\
\Delta_2(c)  & = &
c \otimes S^2a + a \otimes c + b' \otimes Tb - b \otimes Sb
+ Tb \otimes TSb' + b' \otimes TSb' \hspace{.2cm},  \\
\Delta_2(c') & = &
c' \otimes T^2a + a \otimes c' + b \otimes Sb' - b' \otimes Tb'
+ Sb' \otimes STb + b \otimes STb \hspace{.2cm},  \\
\Delta_3(d)  & = &
  d \otimes \ve a + a \otimes d + b \otimes Sc + c' \otimes T^2b'
+ b' \otimes Tc' + c \otimes S^2b \hspace{.2cm},  \\
\Delta_4(a_4) & = & a_4 \otimes T^2a + a \otimes a_4 +
\sum_{g \in P_{24}}
\{s_1(gb) \otimes gSc
+ s_0(ga) \otimes gd
\\ & & 
+ s_2(gc') \otimes gT^2b'
+ s_1(gb') \otimes gTc'
+ s_2(gc) \otimes gS^2b\} \hspace{.2cm}. 
\end{array}$$
\end{Pro}

\begin{Thm} \label{ringP:24}
The ring structure of the group cohomology $H^*(P_{24}; {\mathbb Z})$
is given by 
$H^*(P_{24}; {\mathbb Z}) \iso {\mathbb Z}[\gamma_2,\alpha_4]
/(\gamma_2^2=8\alpha_4, 3\gamma_2=0=24\alpha_4)$.
\end{Thm}

\begin{Thm} \label{ringP:243}
The ring structure of the group cohomology $H^*(P_{24}; {\mathbb Z}_3)$
is given by $H^*(P_{24}; {\mathbb Z}_3) \iso 
{\mathbb Z}_3[\beta_1,\gamma_2,\delta_3,\alpha_4]/(\beta_1^2=0, 
\beta_1\gamma_2=\delta_3, \beta_1\delta_3=0, \gamma_2^2=-\alpha_4, $
$\gamma_2\delta_3=-\beta_1\alpha_4)$.
\end{Thm}

\begin{Thm} \label{ringSMP:24}
Let $M$ be a $3$--dimensional Seifert manifold with $\pi_1(M) \iso P_{24}$.
Then we have the following:\begin{enumerate}
\item
$H^*(M; {\mathbb Z}) \iso
 {\mathbb Z}[\gamma_2, \delta_3]^{\star}/(3\gamma_2 = 0)$.
\item
$H^*(M; {\mathbb Z}_3) \iso
 {\mathbb Z}_3[\beta_1, \gamma_2, \delta_3]^{\star}
 /(\beta_1^2 = 0, \beta_1\gamma_2 = \delta_3)$.
 \end{enumerate}
\end{Thm}

\subsection{Binary octahedral group} \label{P:48}
The $2$-$2$ presentation $P_{48}  = \<T, U\,|\, U^2 = TU^2T, TUT = UTU\>$ 
is given in \cite{cm}.
A more familiar presentation is given by $\< S,T \ | \ S^3 = T^4 = (ST)^2 \>$, 
setting $T = T$,  $U = TS^{-1}$ establishes an isomorphism.

The binary octahedral group $P_{48}$ can be considered as a double
cover of the rotation group of a regular octahedron (or cube), which
is the symmetric group  ${\mathfrak S}_4$.
Thus, there is a commutative diagram of short exact sequences
$$\begin{array}{ccccccccc}
1 & \to & C_2 & \stackrel{\vartriangleleft}{\inc} & P_{48}
              & \stackrel{p}{\onto}               & {\mathfrak S}_4 & \to & 1 \\
 & & \parallel & & \downarrow\subset & & \downarrow\subset & & \\
1 & \to & C_2 & \stackrel{\vartriangleleft}{\inc} & \Spin(3)
              & \stackrel{p}{\onto}               & SO(3)           & \to & 1 
              \hspace{.2cm}.
\end{array}$$
Here, $C_2 = \{1 , z\}$, where $z = T^4 = U^4$, is the centre of $P_{48}$, and
$p(T) = (1\,2\,3\,4)$, $p(U) = (1\,4\,2\,3)$.  One also has
$\left(P_{48}\right)_{ab} \iso \ZZ_2$.

\begin{Pro}\label{P48:cc}
A $4$--periodic resolution ${\mathcal C}$ for $P_{48}$ is given by:
$$\begin{array}{rlrl}
C_0 & =  \<a\>         & \varepsilon(a) &= 1 \hspace{.2cm}, \\
C_1 & =  \<b, b'\> & d_1(b) &= (T-1)a \hspace{.2cm}, \\
        &                           & d_1(b') &= (U-1)a \hspace{.2cm}, \\
C_2 & =  \<c, c'\> & d_2(c) &= (1+TU-U)b + (-1+T-UT)b' \hspace{.2cm}, \\
        &                          & d_2(c') &= (1+TU^2)b + (-1+T-U+TU)b' \hspace{.2cm}, \\
C_3 & =  \<d\>         & d_3(d) &= (1-TU)c + (U-1)c' \hspace{.2cm}, \\
C_4 & =  \<a_4\>   & d_4(a_4) &= Nd \hspace{.2cm}.
\end{array}$$
For any $n \geq 4$, we define $C_n \approx C_{n-4}$ with appropriate subscripts.
\end{Pro}

Let $0 \leq i, j \leq 3$ and $w \in W= \{T^iU^j, T^iUT, T^iU^3T\}_{0
\leq i, j, \leq 3}$.  Then, every word in $P_{48}$ is either in $W$ or
$zW$.  Let $p_i = 1 + T + \cdots + T^{i-1}$ and $q_j = 1 + U + \cdots
+ U^{j-1}$.  In particular, $p_0 = 0 = q_0$.  Write $L$ for $p_4$ and
$M$ for $q_4$.  Define $L' = L(1-U^3)c + (U+p_3TU^3-T^2UT+T^2U^3T)c'$
and $M' = (T+TU+UT+TUT)c + (1-UT)c'$ so that $d_2(L') = (1+z)Lb$ and
$d_2(M') = p_4b - q_4b'$. For further details regarding the above
normal form for the words in $P_{48}$ and the proofs of the formulae
for $d_2(L')$, $d_2(M')$ (which require first deriving further
relations in the group), as well as the proof of the following
proposition (which requires 100 pages of computations) see
\cite{tomoda}.

\begin{Pro}
A contracting homotopy for the chain complex ${\mathcal C}$ above
is given by:
$$\begin{array}{rcl}
s_{-1}(1)      & = & a \hspace{.2cm},  \\
& & \\
s_0(T^iU^ja)   & = & p_ib + T^iq_jb' \hspace{.2cm},  \\
s_0(T^iUTa)    & = & (p_i + T^iU)b + T^ib' \hspace{.2cm},  \\
s_0(T^iU^3Ta)  & = & (p_i + T^iU^3)b + T^iq_3b' \hspace{.2cm},  \\
s_0(za)        & = & (1+TU^2)b + (T+TU+U^2+U^3)b' \hspace{.2cm},  \\
%
s_0(zwa)       & = & s_0(za) + zs_0(wa) \hspace{.2cm},  \mbox{{\rm \ where\ }} w \in W,\\ 
\\
s_1(T^ib)      & = & 0 \hspace{.2cm},  \hfill 0 \leq i \leq 2 \\
s_1(T^3b)      & = & - c' + M' \hspace{.2cm},  \\
s_1(T^iUb)     & = & 0 \hspace{.2cm},  \\
s_1(U^2b)      & = & (U^2-1)c' + (zT^3M'-L') \hspace{.2cm},  \\
s_1(T^iU^2b)   & = & T^{i-1}c' \hspace{.2cm},  \hfill 1 \leq i \leq 3 \\
s_1(T^iU^3b)   & = & 0 \hspace{.2cm},  \\
s_1(UTb)       & = & (U+zT^3U)c + (-1+zT^2+zT^3-zT^3U)c' \\
               &   & + (M'-L') \hspace{.2cm},  \\
s_1(TUTb)      & = & (1+U)c + (-1+U^2)c' + (zT^3M'-L') \hspace{.2cm},  \\
s_1(T^2UTb)    & = & (1+T)c + Uc' \hspace{.2cm},  \\
s_1(T^3UTb)    & = & (T+T^2)c + TUc' \hspace{.2cm},  \\
s_1(U^3Tb)     & = & (U^3+TU^3)c + (-1-TU^3)c' \hspace{.2cm}, \\
s_1(TU^3Tb)    & = & (TU^3+T^2U^3)c + (-TU^2-T^2U^3)c' \hspace{.2cm},  \\
s_1(T^2U^3Tb)  & = & (T^2U^3+T^2U^3T)c - zT^3c' + (L'-T^2M') \hspace{.2cm},  \\
s_1(T^3U^3Tb)  & = & (T^3U^3+T^3U^3T)c - c' + (L'-T^3M') \hspace{.2cm},  \\
\\
s_1(zT^ib)     & = & 0 \hspace{.2cm},  \hfill 0 \leq i \leq 2 \\
s_1(zT^3b)     & = & c' + (L'-M') \hspace{.2cm},  \\
s_1(zT^iUb)    & = & 0 \hspace{.2cm},  \\
s_1(zU^2b)     & = & (T^3+TU^2)c' - TM' \hspace{.2cm},  \\
s_1(zT^iU^2b)  & = & zT^{i-1}c' \hspace{.2cm},  \hfill 1 \leq i \leq 3 \\
s_1(zT^iU^3b)  & = & 0 \hspace{.2cm},  \\
s_1(zUTb)      & = & (T^2+T^3)c + (1+T^2U)c' - M' \hspace{.2cm},  \\
\end{array}$$\def\strutt{\vrule width 0pt height 18pt}
$$\begin{array}{rcl}
s_1(zTUTb)     & = & (T^3+z)c + (1+T^3U)c' - M' \hspace{.2cm},  \\
s_1(zT^2UTb)   & = & (z+zT)c + zUc' \hspace{.2cm},  \\
s_1(zT^3UTb)   & = & (zT+zT^2)c + zTUc' \hspace{.2cm},  \\
s_1(zU^3Tb)    & = & (T^3U^3T + zU^3T)c + (-T+zU^3)c' + (L'-(1+z)M') \hspace{.2cm},  \\
s_1(zTU^3Tb)   & = & (zTU^3 + zTU^3T)c + (1-T^2)c' + (L'-(1+zT)M') \hspace{.2cm},  \\
s_1(zT^2U^3Tb) & = & (zT^2U^3 + zT^3U^3)c - zT^3U^3c' + (L'-zT^2M') \hspace{.2cm},  \\
s_1(zT^3U^3Tb) & = & (U^3 + zT^3U^3)c - U^3c' + (L'-zT^3M') \hspace{.2cm},  \\
\strutt
s_1(T^iU^jb')  & = & 0 \hspace{.2cm},  \hfill 0 \leq i \leq 3, 0 \leq j \leq 2 \\
s_1(T^iU^3b')  & = & -c' + (1-T^i)M' \hspace{.2cm},  \\
s_1(T^iUTb')   & = & -T^ic \hspace{.2cm},  \hfill 0 \leq i \leq 2 \\
s_1(T^3UTb')   & = & -T^3c - c' + M' \hspace{.2cm},  \\
s_1(U^3Tb')    & = & -zT^3U^3c + (U^3-1)c' + (-L'+zT^3M') \hspace{.2cm},  \\
s_1(T^iU^3Tb') & = & -T^{i-1}U^3c + T^iU^3c' + (T^{i-1}-T^i)M' \hspace{.2cm}, 
                   \hfill 1 \leq i \leq 3 \\
\strutt
s_1(zT^iU^jb') & = & 0 \hspace{.2cm},  \hfill 0 \leq i \leq 3, \ 0 \leq j \leq 2 \\
s_1(zT^iU^3b') & = & c' + (L'-(1+zT^i)M') \hspace{.2cm},  \\
s_1(zT^iUTb')  & = & -zT^ic \hspace{.2cm},  \hfill 0 \leq i \leq 2 \\
s_1(zT^3UTb')  & = & -zT^3c + c' + (L'-M') \hspace{.2cm},  \\
s_1(zU^3Tb')   & = & -zU^2c + (1+zU^2)c' + (-1+T^3-z)M' \hspace{.2cm},  \\
s_1(zT^iU^3Tb')& = & -zT^{i-1}U^3c + zT^iU^3c' + z(T^{i-1}-T^i)M' \hspace{.2cm}, 
                   \hfill 1 \leq i \leq 3 \\
\strutt
s_2(T^ic)      & = & 0 \hspace{.2cm},  \\
s_2(Uc)        & = & 0 \hspace{.2cm},  \\
s_2(T^iUc)     & = & -T^{i-1}d \hspace{.2cm},  \hfill 1 \leq i \leq 3 \\
s_2(U^2c)      & = & U^2d \hspace{.2cm},  \\
s_2(T^iU^2c)   & = & T^{i-1}(UT + Tq_3)d \hspace{.2cm},  \hfill 1 \leq i \leq 3 \\
s_2(T^iU^3c)   & = & 0 \hspace{.2cm},  \hfill 0 \leq i \leq 1 \\
s_2(T^2U^3c)   & = & (T + TU + UT)d \hspace{.2cm},  \\
s_2(T^3U^3c)   & = & T^3U^3d \hspace{.2cm},  \\
s_2(UTc)       & = & (1 + U + zT^3 + zT^2UT + zT^3UT)d \hspace{.2cm},  \\
s_2(TUTc)      & = & -Ud \hspace{.2cm},  \\
s_2(T^2UTc)    & = & (-1 + T + UT)d \hspace{.2cm},  \\
s_2(T^3UTc)    & = & (-T + T^2 + TUT)d \hspace{.2cm},  \\
s_2(U^3Tc)     & = & -TU^3d \hspace{.2cm},  \\
s_2(TU^3Tc)    & = & (T + TU + UT - T^2U^3)d \hspace{.2cm},  \\
s_2(T^2U^3Tc)  & = & -(1 + T + U + TU + UT + zT^3UT)d \hspace{.2cm},  \\
s_2(T^3U^3Tc)  & = & -T^3U^3d \hspace{.2cm},  \\
\strutt
s_2(zT^ic)     & = & 0 \hspace{.2cm},  \\
\end{array}$$
$$\begin{array}{rcl}
s_2(zUc)       & = & (T + TU + UT - T^3)d \hspace{.2cm},  \\
s_2(zT^iUc)    & = & -zT^{i-1}d \hspace{.2cm},  \hfill 1 \leq i \leq 2 \\
s_2(zT^3Uc)    & = & zT^3d \hspace{.2cm},  \\
s_2(zU^2c)     & = & (-T - TU - UT + z + zU + T^3UT)d \hspace{.2cm},  \\
s_2(zT^iU^2c)  & = & zT^{i-1}(UT + Tq_3)d \hspace{.2cm},  \hfill 1 \leq i \leq 3 \\
s_2(zU^3c)     & = & z(U^2 + U^3)d \hspace{.2cm},  \\
s_2(zTU^3c)    & = & zp_2U^3d \hspace{.2cm},  \\
s_2(zT^2U^3c)  & = & zT^2U^3d \hspace{.2cm},  \\
s_2(zT^3U^3c)  & = & 0 \hspace{.2cm},  \\
s_2(zT^iUTc)   & = & -T^{i+2}(1 - T - UT)d \hspace{.2cm},  \hfill 0 \leq i \leq 2 \\
s_2(zT^3UTc)   & = & -zT(1 - UT + T^2)d \hspace{.2cm},   \\
s_2(zU^3Tc)    & = & -zU^2d \hspace{.2cm},  
\\
s_2(zTU^3Tc)   & = & -zp_2U^3d \hspace{.2cm},  \\
s_2(zT^2U^3Tc) & = & -(T + TU + UT + zT^2U^3 + zT^3U^3)d \hspace{.2cm},  \\
s_2(zT^3U^3Tc) & = & -U^3d \hspace{.2cm},  \\
\\
s_2(T^ic')     & = & 0 \hspace{.2cm},  \hfill 0 \leq i \leq 2 \\
s_2(T^3c')     & = & -(T + TU + UT)d \hspace{.2cm},  \\
s_2(T^iUc')    & = & 0 \hspace{.2cm},  \\
s_2(U^2c')     & = & 0 \hspace{.2cm},  \\
s_2(T^iU^2c')  & = & T^{i-1}(T + TU + UT)d \hspace{.2cm},  \hfill 1 \leq i \leq 3 \\
s_2(T^iU^3c')  & = & 0 \hspace{.2cm},  \hfill 0 \leq i \leq 2 \\
s_2(T^3U^3c')  & = & (T + TU + UT)d \hspace{.2cm},  \\
s_2(UTc')      & = & z(T^2UT + T^3 + T^3UT)d \hspace{.2cm},  \\
s_2(TUTc')     & = & (-1 - U + UT)d \hspace{.2cm},  \\
s_2(T^2UTc')   & = & (-1 + UT + TUT)d \hspace{.2cm},  \\
s_2(T^3UTc')   & = & (-T + TUT + T^2UT)d \hspace{.2cm},  \\
s_2(U^3Tc')    & = & (zT^3q_3 + U + U^2 + (T^2+T^3)U^3 \\
               &   & + (1+T+zT^2+zT^3)UT + (T+T^2)U^3T)d \hspace{.2cm},  \\
s_2(TU^3Tc')   & = & (1 + p_2U + U^2 + T^3U^3
                     + (2+T+zT^3)UT)d \hspace{.2cm},  \\
s_2(T^2U^3Tc') & = & (-1 + UT + TUT + T^3U^3)d \hspace{.2cm},  \\
s_2(T^3U^3Tc') & = & (-1 - T - TU + TUT - TU^2 - T^3U^3T)d \hspace{.2cm},  \\
\\
s_2(zT^ic')    & = & 0 \hspace{.2cm},  \hfill 0 \leq i \leq 2 \\
s_2(zT^3c')    & = & -(1 + U + zT^3UT)d \hspace{.2cm},  \\
s_2(zT^iUc')   & = & 0 \hspace{.2cm},  \hfill 0 \leq i \leq 1 \\
s_2(zT^2Uc')   & = & z(T^2 + T^3)d \hspace{.2cm},  \\
s_2(zT^3Uc')   & = & (-q_2 + z(T^3 + T^3UT))d \hspace{.2cm},  \\
s_2(zU^2c')    & = & -(T + UT + TU - T^3UT - zq_2)d \hspace{.2cm},  \\
 \end{array}
$$
$$
\begin{array}{rcl}
s_2(zT^iU^2c') & = & z(T^i + T^iU + T^{i-1}UT)d \hspace{.2cm}, 
                     \hfill 1 \leq i \leq 3 \\
s_2(zU^3c')    & = & (T^3U^3 + zU^2)d \hspace{.2cm},  \\
s_2(zTU^3c')   & = & z(U^2 + U^3)d \hspace{.2cm},  \\
s_2(zT^2U^3c') & = & z(U^3 + TU^3)d \hspace{.2cm},  \\
s_2(zT^3U^3c') & = & zT^2U^3d \hspace{.2cm},  \\
s_2(zUTc')     & = & 0 \hspace{.2cm},  \\
s_2(zT^iUTc')  & = & T^i(-T + TUT + T^2UT)d \hspace{.2cm},  \hfill 1 \leq i \leq 2 \\
s_2(zT^3UTc')  & = & z(-T + TUT + T^2UT)d \hspace{.2cm},  \\
s_2(zU^3Tc')   & = & (-1 - T - T^2 - TU - T^2U - zU^2  \\
               &   & - TU^2 - T^2U^2 - T^3U^3T - zU^3T)d \hspace{.2cm},   
\\
s_2(zTU^3Tc')  & = & (-1 + (1+zL)q_3 - zU^2 \\
               &   & + (L+zT^2p_2)U^3 + (p_2+T^3+zL)UT \\
               &   & + (p_3+zT^2p_2)U^3T)d \hspace{.2cm},  \\
s_2(zT^2U^3Tc')& = & (-1 + (1+zTp_3)q_3 + LU^3 \\
               &   & + (p_2+zL)UT + (p_3+zT^3)U^3T)d \hspace{.2cm},  \\
s_2(zT^3U^3Tc')& = & (zT^2p_2 + (1+zT^2p_2)(U+U^2) \\
               &   & + Tp_3U^3 + (p_2+zTp_3)UT + p_3U^3T)d \hspace{.2cm},  \\
\\
s_3(zTU^3Td)   & = & Na_4 \hspace{.2cm},  \\
s_3(gd)        & = & 0 \hspace{.2cm},  \hfill \mbox{ if } g \neq zTU^3T \hspace{.2cm}. 
\end{array}$$
\end{Pro}

\begin{Pro}
A diagonal map $\Delta\co {\mathcal C} \to {\mathcal C}\otimes{\mathcal C}$
for the group $P_{48}$ is given by:
$$\begin{array}{l}
 \Delta_0(a)  = a \otimes a \hspace{.2cm},  \\
 \Delta_1(b)  = b \otimes Ta + a \otimes b \hspace{.2cm},  \\
 \Delta_1(b') = b' \otimes Ua + a \otimes b' \hspace{.2cm},  \\
 \Delta_2(c)  = b \otimes TUb + Tb' \otimes TUb - b' \otimes Ub 
                 + b \otimes Tb' + c \otimes TUTa \hspace{.2cm}  \\
 \hspace{2cm} - Ub \otimes UTb' - b' \otimes UTb' 
                 + a \otimes c \hspace{.2cm},  \\
 \Delta_2(c') = c' \otimes U^2a + b \otimes TU^2b 
                 + Tb' \otimes TU^2b + TUb' \otimes TU^2b 
                 + a \otimes c' \hspace{.2cm}  \\ 
 \hspace{2cm} + b \otimes Tb' - b' \otimes Ub'
                 + b \otimes TUb' + Tb' \otimes TUb' \hspace{.2cm},  \\
 \Delta_3(d)  = a \otimes d - b \otimes TUc + b' \otimes Uc'
                 - Tb' \otimes TUc - c \otimes TUTb' \hspace{.2cm}  \\
 \hspace{2cm} - c \otimes T^2UTb + c' \otimes U^2b' + d \otimes U^3a \hspace{.2cm},  \\
 \Delta_4(a_4)
= \sum_{0 \leq i, j \leq 3} \{p_ib + T^iq_jb'\} \otimes T^iU^jd \\
 
+ \sum_{i = 0}^3 \{(p_i+T^iU)b + T^ib'\} \otimes T^iUTd \\

+ \sum_{i = 0}^3 \{(p_i+T^iU^3)b + T^iq_3b'\} \otimes T^iU^3Td 
+ a \otimes Na_4 \\
- ((z+zT)c + (1+zU)c' + L' - (1+zT^i)M') \otimes c \\

- (zT^2c + zTUc') \otimes Tc \\

- ((U-zT^2+zT^3U)c + (-1+zT^2+zT^3-zT^3U)c' + M' - L') \otimes T^2c  \\

- ((1+U-zT^3)c + U^2c' + (zT^3-1)M') \otimes T^3c \\
\end{array}$$
$$\begin{array}{l}
- ((1+T)c + (U-1)c') \otimes zc  \\

- ((T+T^2)c + (TU-1)c' + (1-T)M') \otimes zTc \\

- ((T^2+T^3)c + T^2Uc' - T^2M') \otimes zT^2c  \\

- ((T^3+z)c + T^3Uc' - T^3M') \otimes zT^3c \\
- (c' + L' - M') \otimes Uc
- (- c' + M') \otimes zUc \\
- ((T^2U^3+T^2U^3T)c - zT^3c' + (L'-T^2M')) \otimes U^2c  \\

- ((T^3U^3+T^3U^3T)c - c' + (L'-T^3M')) \otimes TU^2c \\

- (-zU^2c + (1+zU^2)c' + (-1+T^3-z)M') \otimes T^2U^2c  \\

- (-zU^3c + zTU^3c' + z(1-T)M') \otimes T^3U^2c  
\\

- (-zTU^3c + zT^2U^3c' + z(T-T^2)M') \otimes zU^2c  \\

- (-zT^2U^3c + zT^3U^3c' + z(T^2-T^3)M') \otimes zTU^2c  \\

- ((U^3+TU^3)c + (-1-TU^3)c')    \otimes zT^2U^2c  \\

- ((TU^3+T^2U^3)c + (-TU^2-T^2U^3)c')   \otimes zT^3U^2c  \\
- (c')   \otimes U^3c
- (Tc') \otimes TU^3c
- (T^2c') \otimes T^2U^3c \\

- ((T^3+TU^2)c' - TM') \otimes T^3U^3c \\

- (zc') \otimes zU^3c
- (zTc') \otimes zTU^3c
- (zT2c') \otimes zT^2U^3c
- (U^2b)    \otimes zT^3U^3c \\ 
- (-zT^3c + c' + (L'-M')) \otimes UTc   \\

- (-c) \otimes TUTc
- (-Tc) \otimes T^2UTc
- (-T^2c) \otimes T^3UTc \\

- (-T^3c - c' + M') \otimes zUTc   \\

- (-zc) \otimes zTUTc
- (-zTc) \otimes zT^2UTc
- (-zT^2c) \otimes zT^3UTc \\
- (-U^3c + TU^3c' + (1-T)M') \otimes U^3Tc   \\

- (-TU^3c + T^iU^3c' + (T-T^2)M') \otimes TU^3Tc \\

- (-T^2U^3c + T^3U^3c' + (T^2-T^3)M') \otimes T^2U^3Tc   \\

- (-zU^2c + (1+zU^2)c' + (-1+T^3-z)M') \otimes T^3U^3Tc \\ 

- (-zU^3c + zTU^3c' + z(1-T)M') \otimes zU^3Tc   \\

- (-zTU^3c + zT^2U^3c' + z(T-T^)M') \otimes zTU^3Tc \\

- (-zT^2U^3c + zT^3U^3c' + z(T^2-T^3)M') \otimes zT^2U^3Tc   \\

- (-zT^3U^3c + (U^3-1)c' + (-L'+zT^3M'))    \otimes zT^3U^3Tc \\
+ (c' + (L'-(1+z)M'))   \otimes c' 
+ (c' + (L'-(1+zT)M')) \otimes Tc' \\

+ (c' + (L'-(1+zT^2)M')) \otimes T^2c' 
+ (c' + (L'-(1+zT^3)M')) \otimes T^3c' \\

+ (-c')    \otimes zc' 
+ (-c' + (1-T)M')   \otimes zTc' \\

+ (-c' + (1-T^2)M') \otimes zT^2c' 
+ (-c' + (1-T^3)M') \otimes zT^3c' \\
+ (-zT^3c + c' + (L'-M'))   \otimes UTc' 
+ (-c)     \otimes TUTc' \\

+ (-Tc)       \otimes T^2UTc' 
+ (-T^2c)   \otimes T^3UTc' \\
 
+ (-T^3c - c' + M')     \otimes zUTc' 
+ (-zc)     \otimes zTUTc' \\

+ (-zTc)    \otimes zT^2UTc' 
+ (-zT^2c)  \otimes zT^3UTc' \\
+ (-TU^3c + TU^3c' + (1-T)M')     \otimes U^3Tc' \\
 
+ (-TU^3c + T^2U^3c' + (T-T^2)M')   \otimes TU^3Tc' \\
 
+ (-T^2U^3c + T^3U^3c' + (T^2-T^3)M')   \otimes T^2U^3Tc' \\
 \end{array}$$
$$\begin{array}{l}
+ (-zU^2c + (1+zU^2)c' + (-1+T^3-z)M')  \otimes T^3U^3Tc'\\
 
+ (-zU^3c + zTU^3c' + z(1-T)M')     \otimes zU^3Tc' \\
 
+ (-zTU^3c + zT^2U^3c' + z(T-T^2)M')    \otimes zTU^3Tc' \\

+ (-zT^2U^3c + zT^3U^3c' + z(T^2-T^3)M')    \otimes zT^2U^3Tc' \\
 
+ (-zT^3U^3c + (U^3-1)c' + (-L'+zT^3M'))      \otimes zT^3U^3Tc' \\
+ \{ -(T + UT + TU - T^3 -zU)d \} \otimes b' 
+ \{ z(1 + TU)d \} \otimes Tb' \\

+ \{ z(T + T^2 + T^3 + T^2U)d \} \otimes T^2b'  
\\

+ \{ (1 + U - zT^3U + zT^3UT)d \} \otimes T^3b' \\

+ \{ Ud \} \otimes zb'
+ \{ (1 + TU)d \} \otimes zTb' 
+ \{ (T + T^2U)d \} \otimes zT^2b' \\

+ \{ (T^2 + T^3U)d \} \otimes zT^3b' \\
+ \{ zU^2d \} \otimes Ub' 
%
%
- \{ (T + TU + UT)d \} \otimes zT^3Ub' \\
+ \{ T^3U^3d \} \otimes U^2b' 
+ \{ zU^2d \} \otimes TU^2b' 
+ \{ zp_2U^3d \} \otimes T^2U^2b'  \\

+ \{ z(T^2U^3 + T^3U^3)d \} \otimes T^3U^2b' 
+ \{ U^3d \} \otimes zU^2b' 
+ \{ TU^3d \} \otimes zTU^2b' \\

- \{ (T + TU + UT - T^2U^3)d \} \otimes zT^2U^2b' 
+ \{ (T + TU + UT)d \} \otimes zT^3U^2b'  \\
+ \{ d \} \otimes U^3b' 
+ \{ Td \} \otimes TU^3b' 
+ \{ T^2d \} \otimes T^2U^3b' \\

- \{ (T + TU + UT - T^3)d \} \otimes T^3U^3b' 
+ \{ zd \} \otimes zU^3b' 
+ \{ zTd \} \otimes zTU^3b'\\
 
+ \{ zT^2d \} \otimes zT^2U^3b' 
+ \{ (-q_2 + z(T^3 + T^3UT))d \} \otimes zT^3U^3b' \\
+ \{ (-T^3 + T^3UT + zUT)d \} \otimes UTb' 
+ \{ z(-T + TUT + T^2UT)d \} \otimes TUTb' \\

+ \{ z(T^2UT + T^3 + T^3UT)d \} \otimes T^2UTb' 
+ \{ (-1 - U + UT)d)   \otimes T^3UTb' \\

+ \{ (-1 + UT + TUT)d \} \otimes zUTb' 
+ \{ (-T + TUT + T^2UT)d \} \otimes zTUTb' \\

+ \{ (-T^2 + T^2UT + T^3UT)d \} \otimes zT^3UTb' \\
+ \{ -(1 + TUT + Tq_3 + T^3U^3)d \} \otimes U^3Tb' \\

+ \{ -(1 + (T+T^2)q_3 + T^3U^3T)d \} \otimes TU^3Tb' \\

+ \{ -(1 + (T+T^2+T^3)q_3 + zU^2 + T^2UT + T^3U^3T + zU^3T)d \}
      \otimes T^2U^3Tb' \\
 
+ \{ (z(T+T^2+T^3) + (1+zT+zT^2+zT^3)U + (1+zT+zT^2+zT^3)U^2 \\ 
 \hspace{1cm} + (L+zT^2+zT^3)U^3 + (p_2+T^3+zL)UT + p_3U^3T)d \} 
      \otimes T^3U^3Tb' \\

+ \{ (z(T^2+T^3) + (1+zT^2+zT^3)U + (1+zT^2+zT^3)U^2 + LU^3 \\
 \hspace{1cm} + (p_2+zT+zT^2+zT^3)UT + (p_3+zT^3)U^3T)d \}
      \otimes zU^3Tb' \\

+ \{ (zT^3 + (1+zT^3)U + (1+zT^3)U^2 + (T+T^2+T^3)U^3 \\
 \hspace{1cm} + (1+T+zT^2+zT^3)UT + p_3U^3T)d \} 
      \otimes zTU^3Tb' \\

+ \{ (U + U^2 + (T^2+T^3)U^3 + (1+T+zT^3)UT \\
 \hspace{1cm} + (T+T^2)U^3T)d \}
      \otimes zT^2U^3Tb' \\

+ \{ (1 + (1+T)U + T^3U^3 + (2+T+zT^3)UT)d \} \otimes zT^3U^3Tb' \\
+ \{ zd \} \otimes b 
+ \{ zTd \} \otimes Tb 
+ \{ zT^2d \} \otimes T^2b 
+ \{ d \} \otimes zb \\

+ \{ Td \} \otimes zTb 
+ \{ T2d \} \otimes zT^2b 
- \{ (T + TU + UT - T^3)d \} \otimes zT^3b \\
+ \{ (T^3 - z - T^3UT)d \} \otimes Ub 
+ \{ z(1 - T - UT)d \} \otimes TUb \\
 
+ \{ z(T - TUT + T^3)d \} \otimes T^2Ub \\
 
- \{ (1 + U + zT^3 + zT^2UT + zT^3UT)d \} \otimes T^3Ub \\
 
+ \{ Ud \} \otimes zUb 
- \{ (-1 + T + UT)d \} \otimes zTUb \\
 \end{array}$$
$$\begin{array}{l}
- \{ (-T + T^2 + TUT)d \} \otimes zT^2Ub 
+ \{ (T^2 - T^3 - T^2UT)d \} \otimes zT^3Ub \\
- \{ T^3U^3d \} \otimes U^2b 
- \{ z(U^2 + U^3)d \} \otimes TU^2b 
- \{ zp_2U^3d \} \otimes T^2U^2b \\

- \{ zT^2U^3d \} \otimes T^3U^2b 
- \{ (T + TU + UT)d \} \otimes zT^3U^2b  
\\
+ \{ T^3U^3d \} \otimes U^3b 
+ \{ zU^2d \} \otimes TU^3b 
+ \{ zp_2U^3d \} \otimes T^2U^3b \\

+ \{ (T + TU + UT + zT^2U^3 + zT^3U^3)d \} \otimes T^3U^3b 
+ \{ U^3d \} \otimes zU^3b \\
 
+ \{ TU^3d \} \otimes zTU^3b 
- \{ (T + TU + UT - T^2U^3)d \} \otimes zT^2U^3b \\

+ \{ (1 + T + U + TU + UT + zT^3UT)d \} \otimes zT^3U^3b  \\
- \{ (TUT + T^2q_3)d \} \otimes U^3Tb 
- \{ (T^2UT + T^3q_3)d \} \otimes TU^3Tb \\
- \{ (-T - TU - UT + z + zU + T^3UT)d \} \otimes T^2U^3Tb \\
- \{ z(UT + Tq_3)d \} \otimes T^3U^3Tb 
- \{ z(TUT + T^2q_3)d \} \otimes zU^3Tb \\
- \{ z(T^2UT + T^3q_3)d \} \otimes zTU^3Tb 
- \{ U^2d \} \otimes zT^2U^3Tb \\
- \{ (UT + Tq_3)d \} \otimes zT^3U^3Tb \\
+ Na_4 \otimes T^2U^3Ta \hspace{.2cm}.  \\
\end{array}$$
\end{Pro}

\begin{Thm} \label{ringP:48}
The ring structure of the group cohomology $H^*(P_{48}; {\mathbb Z})$
is given by 
$H^*(P_{48}; {\mathbb Z}) \iso {\mathbb Z}[\gamma_2, \alpha_4]
/(\gamma_2^2 = 24\alpha_4, 2\gamma_2 = 0 = 48\alpha_4)$.
\end{Thm}

As mentioned in the Preliminaries (\fullref{sec:2}), we will give the proof of the next theorem
to provide an example of how the cohomology ring is determined from the
resolution and diagonal map.

\begin{Thm} \label{ringP:482}
The ring structure of the group cohomology $H^*(P_{48}; {\mathbb Z}_2)$
is given by 
$H^*(P_{48}; {\mathbb Z}_2) \iso
{\mathbb Z}_2[\beta_1, \gamma_2, \delta_3, \alpha_4]
/(\beta_1^2 = \gamma_2, 
\gamma_2^2 = 0 = \delta_3^2, \beta_1\gamma_2 = \delta_3, 
\beta_1\delta_3 = 0 = \gamma_2\delta_3)$.
\end{Thm}
\Proof
We consider the coefficients $\ZZ_2$ as an $R$--algebra with trivial $R$--module
structure.  The cochain complex $\hom_R({\mathcal C}, \ZZ_2)$ is then generated
by the dual classes $\hat{a}$, $\hat{b}$, $\hat{b'}$, $\hat{c}$, 
$\hat{c'}$, $\hat{d}$, and $\hat{a_4}$, where, for example,  
$\hat{b}(b) = 1$, $\hat{b}(b') = 0$, etc.  We find, 
for the coboundary $\partial$, 
\begin{align*}\left(\partial\hat{a}\right)(b) &= \hat{a}(d_1b) = \hat{a}(Ta+a) = 1+1 = 0\\
\left(\partial\hat{a}\right)(b') &= \hat{a}(d_1b') = \hat{a}(Ua+a) = 1+1 
= 0 \ \ \ ,\tag*{\hbox{and}}
\\ \partial\hat{a} &= 0. \tag*{\hbox{hence,}}
\\
\partial\hat{b} & =  \hat{c}  ,\qquad  \partial\hat{b'} =  \hat{c}
 \hspace{.2cm}, \tag*{\hbox{Similarly,}}\\
\partial\hat{c} & = 0, \qquad \partial\hat{c'}  =  0
 \hspace{.2cm}, \\
\partial\hat{d} & =  0, \\
\partial\hat{a_4} & =  0. 
\end{align*}
The cohomology therefore has generating classes and representative cocycles as
shown in the following table.
$$\begin{array}{|c|l|}
\hline
{\rm Dimension} & {\rm Cohomology\ class\ \&\ representative\ cocycle} \st \\
\hline
0 & 1 = [\hat{a}] \st \\
1 & \beta_1 = [\hat{b} + \hat{b'}] \\
2 & \gamma_2 = [\hat{c'}] \\
3 & \delta_3 = [\hat{d}] \\
4 & \alpha_4 = [\hat{a_4}] = {\rm periodicity\ class} \\
\hline
\end{array}$$
Since $(\hat{b}+\hat{b'})\otimes(\hat{b}+\hat{b'})
= \hat{b}\otimes\hat{b} + \hat{b}\otimes\hat{b'}
+ \hat{b'}\otimes\hat{b} +  \hat{b'}\otimes\hat{b'}$,
it follows that $\beta_1^2 = \lambda\gamma_2$, where $\lambda$ is the
number of terms (mod $2$) in $\Delta(c')$ of the form $xb\otimes yb$, 
$xb\otimes yb'$, $xb'\otimes yb$, and $xb'\otimes yb'$, for any 
$x, y, \in P_{24}$.  Using \fullref{P48:cc}, a simple count shows
that $\lambda = 7 = 1$.  Thus, $\beta_1^2 = \gamma_2$.  The cup products
$\gamma_2^2 = 0$, $\beta_1\gamma_2 = \delta_3$, and 
$\beta_1\delta_3 = 0$ are computed similarly.  Then, 
$\gamma_2\delta_3 = \beta_1^2\delta_3 = \beta_1(\beta_1\delta_3) = 0$
as well as $\delta_3^2 = \beta_1^2\gamma_2^2 = 0$.  
Periodicity then determines all further cup products.
\end{proof}

\begin{Thm} \label{ringP:483}
The ring structure of the group cohomology $H^*(P_{48}; {\mathbb Z}_3)$
is given by 
$H^*(P_{48}; {\mathbb Z}_3) \iso
{\mathbb Z}_3[\delta_3, \alpha_4]/(\delta_3^2 = 0)$. For $p > 3$, 
$H^*(P_{48};\ZZ_p) \iso \ZZ_p[\alpha_4]$.
\end{Thm}

\begin{Thm} \label{ringSMP:48}
Let $M$ be a $3$--dimensional Seifert manifold with $\pi_1(M) \iso P_{48}$.
Then we have the following:\begin{enumerate}
\item
$H^*(M; {\mathbb Z}) \iso {\mathbb Z}[\gamma_2, \delta_3]^{\star}/(2\gamma_2 = 0)$. 
\item
$H^*(M; {\mathbb Z}_2) \iso
{\mathbb Z}_2[\beta_1, \gamma_2, \delta_3]^{\star}
/(\beta_1^2 = \gamma_2, \beta_1\gamma_2 = \delta_3)$.
\item
$H^*(M; {\mathbb Z}_p) \iso {\mathbb Z}_p[\delta_3]^{\star}$,
for $p \neq 2$.
\end{enumerate}
\end{Thm}

\subsection{Binary icosahedral group} \label{P:120}
Following Coxeter--Moser \cite{cm}, the presentation we use for the binary icosahedral group is
$P_{120} = \<A, B\,|\, AB^2A=BAB, BA^2B=ABA \>$.  This is the
fundamental group of the homology sphere discovered by Poincar\'e,
and this is the only known homology $3$--sphere with a finite
fundamental group.  Of course, the fact that $H_1(P_{120}; \ZZ) = 0$
(and hence it is a homology sphere) follows from
$\left(P_{120}\right)_{ab} = 0$.  Once again, it can be regarded as a
double cover, in this case, of the simple group ${\mathfrak A}_5$
(which is the rotation group of a regular icosahedron or dodecahedron),
as shown by the commutative diagram
$$\begin{array}{ccccccccc}
1 & \to & C_2 & \stackrel{\vartriangleleft}{\inc} & P_{120}
              & \stackrel{p}{\onto}               & {\mathfrak A}_5 & \to & 1 \\
 & & \parallel & & \downarrow\subset & & \downarrow\subset & & \\
1 & \to & C_2 & \stackrel{\vartriangleleft}{\inc} & \Spin(3)
              & \stackrel{p}{\onto}               & SO(3)           & \to & 1 
              \hspace{.2cm}.
\end{array}$$
Here, we can take $p(A) = (1\,2\,3\,4\,5)$,
$p(B) = (1\, 3\,4\,2\,5)$, and $C_2 = \{1, z\}$ where 
$z := (ABA)^3 = (BAB)^3$.

\begin{Pro}
A $4$--periodic resolution ${\mathcal C}$ for $P_{120}$ is given by:
$$\begin{array}{rlrl}
C_0 & =  \<a\>         & \varepsilon(a) &= 1 \hspace{.2cm}, \\
C_1 & =  \<b, b'\> & d_1(b) &= (A-1)a \hspace{.2cm}, \\
               &                     & d_1(b') &= (B-1)a \hspace{.2cm}, \\
C_2 & =  \<c, c'\> & d_2(c) &= (1-B+AB^2)b + (-1+A+AB-BA)b' \hspace{.2cm}, \\
               &                     & d_2(c') &= (-1+B+BA-AB)b + (1-A+BA^2)b' \hspace{.2cm}, \\
C_3 & =  \<d\>         & d_3(d) &= (1-BA)c + (1-AB)c' \hspace{.2cm}, \\
C_4 & =  \<a_4\>   & d_4(a_4) &= Nd \hspace{.2cm}.
\end{array}$$
For any $n \geq 4$, we define $C_n \approx C_{n-4}$ with appropriate subscripts.
\end{Pro}

For this group, the construction of a contracting homotopy $s$ seems
daunting, since the corresponding work for $P_{48}$ took nearly 100 pages.
However, exactness of the resolution ${\mathcal C}$ has been verified using
a computer (cf \fullref{sec:2}).

\begin{Pro}
A diagonal map $\Delta\co {\mathcal C} \to {\mathcal C}\otimes{\mathcal C}$, 
through dimension $3$,
for the group $P_{120}$ is given by:
\begin{eqnarray*}
& & \Delta_0(a)  = a \otimes a \hspace{.2cm},  \\
& & \Delta_1(b)  = b \otimes Aa + a \otimes b \hspace{.2cm},  \\
& & \Delta_1(b') = b' \otimes Ba + a \otimes b' \hspace{.2cm},  \\
& & \Delta_2(c)  = a \otimes c + c \otimes BABa - b' \otimes Bb
                 - b' \otimes BAb' - Bb \otimes BAb'  \\
& & \hspace{1cm} + ABb' \otimes AB^2b + b \otimes Ab' + Ab' \otimes ABb' 
                 + Ab' \otimes AB^2b \\
& & \hspace{1cm} + b \otimes AB^2b + b \otimes ABb' \hspace{.2cm},  \\
& & \Delta_2(c') = a \otimes c' + c' \otimes ABAa - b \otimes Ab'
                 - b \otimes ABb - Ab' \otimes ABb \\ 
& & \hspace{1cm} + BAb \otimes BA^2b'  + b' \otimes Bb + Bb \otimes BAb
                 + Bb \otimes BA^2b'\\
& & \hspace{1cm} + b' \otimes BA^2b' + b' \otimes BAb \hspace{.2cm},  \\
& & \Delta_3(d)  = a \otimes d + d \otimes (BA)^2Ba - c \otimes BABb
                 - b \otimes ABc' - Ab' \otimes ABc' \\
& & \hspace{1cm} - Bb \otimes BAc - b' \otimes BAc - c' \otimes (AB)^2b 
                 - c' \otimes ABAb' \\
& & \hspace{1cm} - c \otimes (BA)^2b' \hspace{.2cm}. 
\end{eqnarray*}
\end{Pro}

We remark that for computing ring structures of $H^*(P_{120}; A)$ with
twisted coefficients $A$, one probably requires an explicit formulation
of $\Delta_4$.

From the above, we have the next theorem.
\begin{Thm} \label{ringP:120}
\begin{eqnarray*}
H^l(P_{120}; {\mathbb Z}) & = &
\left\{\begin{array}{ll}
{\mathbb Z} \hspace{.2cm},        & \mbox{ if } l = 0, \\
0 \hspace{.2cm},                  & \mbox{ if } l \not\equiv 0 (\mbox{ mod}\, 4), \\
{\mathbb Z}_{120} \hspace{.2cm},  & \mbox{ if } l \equiv 0 (\mbox{ mod}\, 4)
                    \mbox{ and } l > 0.
\end{array}\right.
\end{eqnarray*}
It follows that
$H^*(P_{120}; {\mathbb Z}) \iso {\mathbb Z}[\alpha_4]/(120\alpha_4 = 0)$ 
and $H^*(P_{120}; \ZZ_n) \iso \ZZ_n[\alpha_4]$  with $n$  any divisor of $120$.
Also $H^*(M; \ZZ) = \ZZ[\delta_3]^{\star}$,   where
$M = S^3/P_{120}$ is the Poincar\'e homology sphere.
\end{Thm}

\section[The groups P'(8.3^k) and B(2^k(2n+1))]{The groups $P'_{8\cdot3^k}$ and $B_{2^k(2n+1)}$} \label{sec:5}

In this chapter we compute the ring structures of cohomology groups of the
groups $P'_{8\cdot3^k}$ and $B_{2^k(2n+1)}$.
For these groups, we employ a more geometrical approach, using appropriate 
Seifert manifolds.  We assume some general familiarity with Seifert manifolds
(good references are Seifert \cite{seifert}, Hempel \cite{hempel} and Orlik \cite{orlik}), and will
merely introduce Seifert's notation for them.  One writes
$$
M = \left(\{O, N\}, \{o, n\}, g: e; (a_1, b_1), \cdots, (a_q, b_q)\right),
$$
where \fullref{tblSeifert} describes the meaning of each symbol.
\begin{table}[h] \label{tblSeifert} \centering 
\begin{tabular}{|rl|}
\hline
\st
$\{O, N\}$:  & the orientability of the Seifert manifold $M$:       \\
\st
             & $O$ means that $M$ is orientable, and                \\
\st
             & $N$ means that $M$ is nonorientable,                \\
\st
$\{o, n\}$:
\st
             & the orientability of its orbit surface $V$:          \\
\st
             & $o$ means that $V$ is orientable, and                \\
\st
             & $n$ means that $V$ is nonorientable                 \\
\st
$g$:         & if $o$, then $g \geq 0$ equals the genus of $V$,            \\
\st
             & if $n$, then $g \geq 1$ equals number of cross-caps of $V$, \\
\st
$e$:         & the Euler number, obtained from a regular fibre;     \\
\st
$q$:         & the number of singular fibres;                    \\

\st
$(a_i, b_i)$:
\st
             & the relatively prime integer pairs characterizing the $i$--th \\
\st
             & singular fibre with $0 < b_i < a_i$.           \\ 
\hline
\end{tabular}
\vskip8pt
\caption{Presentation of Seifert manifolds} 
\end{table} 
The resolutions of ${\mathbb Z}$ over the group ring $R$, and the
diagonal $\Delta$, are 
based on the methods of Bryden, Hayat-Legrand, Zieschang and Zvengrowski
in \cite{bhzz1,bhzz2,bz}, appropriately modified to account for 
the universal cover $\wwtilde{M}$ now being $S^3$ instead of $\RR^3$.
We verified that the chain complexes provided below
are indeed resolutions by using a computer program in GAP, at least for small orders.

In \fullref{tblSMfundgps}, we list Seifert manifolds with finite fundamental group and the
corresponding presentation as a Seifert manifold. This table is based on 
Orlik \cite[p\,112]{orlik}, with minor notational changes and three small corrections:
in the first case $n\not= 0$ is added, otherwise $\pi_1(M)\approx {\mathbb Z}$ is infinite,
in the second case $B_{2^{k+3}a_3}$ is incorrectly given as  $B_{2^{k+2}a_3}$ in \cite{orlik},
and in the third case the equation $m = 3^{k-1}m'$ is incorrectly given as $m = 3^km'$ in \cite{orlik}.
Finiteness implies $g = 0$ in the $(O, o)$ case and $g = 1$ in the $(O, n)$ case. Denoting
the number of singular fibres by $q$,
 the fundamental groups \cite[\S2]{bz} are then given by
\begin{align*} \pi_1(M) &= \< \  s_1,...,s_q,h \ |\ [s_j,h], s_j^{a_j}h^{b_j},  
s_1\cdots s_q h^{-e} \ \> , \ (O,o)-{\rm case},  \\
 \pi_1(M) &= \< \ s_1,...,s_q,h,v \ |\ [s_j,h], s_j^{a_j}h^{b_j}, vhv^{-1}h, 
s_1\cdots s_q v^2  h^{-e} \ \> , \ (O,n)-{\rm case} \ . \end{align*} 
Note that the same group may appear more than once since fibre-inequivalent Seifert
spaces can have the same fundamental group, this is characteristic of ``small" Seifert manifolds
\cite[p\,91]{orlik}. Also note $B_{2^k(2n+1)}$  is defined for  $n\geq 0$, with the 
group isomorphic to  ${\mathbb Z}_{2^k}$  when  $n = 0$, and for  $k=2$ there is an
isomorphism  $B_{4\cdot (2n+1)} \iso Q_{4\cdot (2n+1)}$ \cite{milnor}.  For this 
reason, in \fullref{bgroups} below, we only consider
 $B_{2^k(2n+1)}$  for  $n \geq 1$, $k \geq 3$.

\begin{table}[h] \label{tblSMfundgps} \centering
\begin{tabular}{|ll|}%
\hline\st
Seifert structure & Fundamental Group \\
\hline\st
$(O, o, 0: e; (a_1, b_1), (a_2, b_2))$
  & $\pi_1(M)\approx  {\mathbb  Z}_n$, \ \
    $n: = |e a_1 a_2 + a_1 b_2 + b_1 a_2|,$  \\
  & $n \not= 0$ ($(a_j,b_j) = (1 , 0)$ is allowed). \\
$(O, o, 0: e; (2, 1), (2, 1), (a_3, b_3))$
  & Let $m = |(e+1)a_3 + b_3|$. \\
  & If $m$ is odd, then $\pi_1(M)\approx {\mathbb  Z}_m \times Q_{4a_3}$. \\
  & If $m$ is even, then $4|m$ and $(a_3, 2) = 1$. \\
  & Set $m = 2^{k+1}m''$, $m''$ odd. \\   
  & Then  $\pi_1(M)\approx {\mathbb  Z}_{m''} \times B_{2^{k+3}a_3}$. \\
$(O, o, 0: e; (2, 1), (3, b_2), (3, b_3))$
  & Let $m = |6e + 3 + 2(b_2+b_3)|$ and \\
  & $m = 3^{k-1}m'$ with $(m', 3) = 1$. \\
  & If $k = 1$, then $(m, 6) = 1$, $b_2 = 1 = b_3$,\\
  & and $\pi_1(M)\iso {\mathbb  Z}_m \times P'_{24} \iso {\mathbb  Z}_m \times P_{24}$. \\
  & If $k \geq 2$, then $(m', 6) = 1$, $b_2 = 1$, $b_3 = 2$, \\
  & and $\pi_1(M)\iso {\mathbb  Z}_{m'} \times P_{8\cdot3^k}'$. \\
$(O, o, 0: e; (2, 1), (3, b_2), (4, b_3))$
  & $\pi_1(M)\iso{\mathbb  Z}_m \times P_{48}$, \\ 
  & where $m = |12e + 6 + 4b_2 + 3b_3|$. \\
$(O, o, 0: e; (2, 1), (3, b_2), (5, b_3))$
  & $\pi_1(M)\iso{\mathbb  Z}_m \times P_{120}$, \\  
  & where $m = |30e + 15 + 10b_2 + 6b_3|$. \\
$(O, n, 1: e; (a_1, b_1))$
  & Let $m = |e a_1 + b_1|$. \\
  & If $a_1$ is odd, then $\pi_1(M)\iso {\mathbb  Z}_{\alpha_1} \times Q_{4m}$. \\
  & If $a_1$ is even, then
    $\pi_1(M)\iso {\mathbb  Z}_{a_1'} \times B_{2^{k+2}m}$, \\
  & where $a_1 = 2^k a_1'$, $k \geq 1$, and $(a_1', 2) = 1$. \\
\hline
\end{tabular}  
\vskip8pt
\caption{Seifert manifolds with finite fundamental groups 
(following Orlik \protect\cite[p\,112]{orlik})}
\end{table}

The following elements in $G = \pi_1(M)$, $R = \ZZ G$, and in ${\mathcal C}$ 
 are necessary to define the resolutions for these groups used in \fullref{pgroups} and \fullref{bgroups}.
\begin{enumerate}
\item (in $G$, $M = (O,o,0 : e ; (a_1, b_1), ... (a_q, b_q))$)\qua
Choose positive integers $c_j$ and $d_j$
satisfying $a_jd_j - b_jc_j = 1$, $1 \leq j \leq q $, and let
$t_j = \smash{s_j^{\smash{c_j}}h^{\smash{d_j}}}$. Also define $a_0 = 1$, $b_0 = e$, $c_0 = 1$,  $d_0 = e+1$, \
and $s_0 = h^{-e}$.  As a consequence, $a_0d_0 - b_0c_0 = 1$, $t_0 = \smash{s_0^{\smash{c_0}}h^{\smash{d_0}}}
= h$, $s_j = \smash{t_j^{\smash{-b_j}}}$,  $\smash{h = t_j^{\smash{a_j}}}$,  $0 \leq j \leq q$.

\item (in $G$, $M = (O,n,1: e,(a_1, b_1), ... (a_q, b_q)) $)\qua 
The relation  $vhv^{-1}h = 1$  implies  $h^i v h^i = v$, $i \in \ZZ,$ as well as 
 $h v^2 = v^2 h$.

\item (in $G$)\qua  Let  $r_{-1} = 1$, $r_j = s_0s_1\cdots s_j,   0 \leq j \leq q, $
in the $(O,n)$--case also $r_{q+1} =  s_0s_1\cdots s_q v^2$.

\item (in ${\mathcal C}$)\qua
Let $\pi_j^1 = r_{j-1}(\sigma_j^1+\rho_j^1) - r_j\sigma_j^1$, 
in the $(O, n)$--case, $\pi_{q+1}^1 = r_{q}(1+v)\nu_1^1$.

\item (in ${\mathcal C}$)\qua
Let $\pi_j^2 = -r_{j-1}(\sigma_j^2+\rho_j^2) + r_j\sigma_j^2$, 
in the $(O, n)$--case, $\pi_{q+1}^2 = r_{q}(hv-1)\nu_1^2$.

\item (in ${\mathbb Z}G$)\qua
Let $F_j = \upnpfrac{t_j^{a_j}-1}{t_j-1}$ and
$G_j = \upnpfrac{1-t_j^{-b_j}}{t_j-1}$.
\end{enumerate}

\subsection[The groups P'(8.3^k)]{The groups $P'_{8\cdot3^k}$} \label{pgroups}

The group $P'_{8\cdot3^k}$, $k \geq 1$, are given by the following
presentation:
\begin{eqnarray*}
P'_{8\cdot3^k}
 & = & \< x, y, z\,|\, x^2=(xy)^2=y^2, zxz^{-1}=y, zyz^{-1}=xy, z^{3^k}=1 \>. 
\end{eqnarray*}
One can also represent these groups as semidirect products; namely,
\linebreak
$P'_{8\cdot3^k}  \iso Q_8 \rtimes C_{3^k}$.  Equivalently, one has a
split short exact sequence
\begin{center}\begin{picture}(150, 10)
\thinlines
\put(  5,0){\makebox(0,0)[t]{\smash{$1$}}}
\put( 10,3){\vector(1,0){10}}
\put( 30,0){\makebox(0,0)[t]{\smash{$Q_8$}}}
\put( 44,-2){\makebox(0,0)[t]{\smash{\Large{$\hookrightarrow$}}}}
\put( 43,5){\makebox(0,0)[t]{\smash{\small{$\vartriangleleft$}}}}
\put( 65,0){\makebox(0,0)[t]{\smash{$P'_{8\cdot3^k}$}}}
\put( 84,2){\makebox(0,0)[t]{\smash{\Large{$\twoheadrightarrow$}}}}
\put( 84,10){\makebox(0,0)[t]{\smash{\small{$p$}}}}
\qbezier( 79,2)(86,0)(93,2)\put(76,3){\vector(-3,1){1}}
\put( 84,-6){\makebox(0,0)[t]{\smash{\small{$s$}}}}
\put(105,0){\makebox(0,0)[t]{\smash{$C_{3^k}$}}}
\put(115,3){\vector(1,0){10}}
\put(130,0){\makebox(0,0)[t]{\smash{$1$}}}
\end{picture}\end{center}
where $Q_8$ is the subgroup generated by $x$, $y$, $C_{3^k}$ is the
cyclic group with $z$ as generator, $p(x, y) = 1$, $p(z) = z$, and
$s(z) = z$.  We remark that $P'_{8\cdot3} \iso P_{24}$, the binary
tetrahedral group, and $(P'_{8\cdot3^k})_{ab} = \ZZ_{3^k}$.

Since $P'_{8\cdot3} \iso P_{24}$, we are only concerned with the case $k \geq 2$
(also, as shown by \fullref{tblSMfundgps}, the Seifert structure is slightly different when $k = 1$).
Let $M = (O, o, 0: e; (2, 1), (3, 1), (3, 2))$ where 
$e = (1/2)(3^{k-2} - 3)$, $k\geq2$. Then, again following \mbox{\fullref{tblSMfundgps}},
$m' = 1$  and  $\pi_1(M) \iso  P'_{8\cdot3^k}$. We now outline a proof of this, partly 
because none of the isomorphisms in \fullref{tblSMfundgps} are explicitly proved in \cite{orlik}, and also
because \cite{orlik} has a minor error in this case.

\begin{Pro}   
With $M$ and $e$ as above,  $\pi_1(M) \iso  P'_{8\cdot3^k}$.
\end{Pro}

\begin{proof}[Proof outline]
 The fundamental group $\pi_1(M)$ is given by
$$\pi_1(M) = \<s, t, u, h | [s,h] = [t,h] = [u,h] = s^2h = t^3h = u^3h^2 = stuh^{-e} = 1\>.$$
Here, we have used $s$, $t$, $u$ instead of the notation $s_1$, $s_2$, $s_3$ used
in \cite{bhzz2}.  
Let $m = 3^{k-1}$ and $n = 3e+5$, and define $\varphi\co\pi_1(M) \to P'_{8\cdot3^k}$ by
$$
\varphi(s) = x^{2e+1}z^{3(7-3e^2)},\ \ \ 
\varphi(t) = x^3z^n,\ \ \ 
\varphi(u) = z,\ \ \ 
\varphi(h) = x^2z^{3(n-1)}
$$
and $\psi\co P'_{8\cdot3^k}\to \pi_1(M)$ by
$$
\psi(x) = s^m,\ \ \ 
\psi(y) = s^{2m-3}tst^2,\ \ \ 
\psi(z) = u.
$$
One can show (for full details see Tomoda \cite{tomoda}) that the maps $\varphi$ and $\psi$ are 
well defined and inverse 
isomorphisms between 
the fundamental group $\pi_1(M)$ and the group $P'_{8\cdot3^k}$, $k\geq2$. \end{proof}

\begin{Pro}\label{prop5.2}
A resolution $\mathcal{C}$ for $P'_{8\cdot3^k}$ is given,
with  $0\leq j\leq3$, by:
$$\begin{array}{lcl}
C_0 & = & \<\sigma_j^0\> \hspace{.2cm}, \\
C_1 & = & \<\sigma_j^1, \rho_j^1, \eta_j^1\> \hspace{.2cm}, 
\mbox{ with } \sigma_0^1 = 0 \\
C_2 & = & \<\sigma_j^2, \rho_j^2, \mu_j^2, \delta^2\> \hspace{.2cm}, 
\mbox{ with } \sigma_0^2 = 0 \\
C_3 & = & \<\sigma_j^3, \delta^3\> \hspace{.2cm}, \\
C_4 & = & \<\sigma_0^4\>, 
\end{array}$$
along with
$$\begin{array}{rlrl}
d_1(\sigma_j^1)     & =  \sigma_j^0-\sigma_0^0 \hspace{.2cm}, 
 & d_1(\rho_j^1) & =  (s_j-1)\sigma_0^0  \hspace{.2cm}, \\
d_1(\eta_j^1)         & =  (h-1)\sigma_0^0 \hspace{.2cm}, 
&   &      \\
d_2(\sigma_j^2)     & =  \eta_0^1 - \eta_j^1+ (h-1)\sigma_j^1 \hspace{.2cm}, 
& d_2(\rho_j^2) & =  (1-s_j)\eta_j^1 + (h-1)\rho_j^1 \hspace{.2cm},  \\
d_2(\mu_j^2)        & =  F_j\rho_j^1 + G_j\eta_j^1 \hspace{.2cm}, 
   & d_2(\delta^2) & =  \sum_{j=0}^3\pi_j^1 \hspace{.2cm},  \\
d_3(\sigma_j^3)     & =  \rho_j^2 + (1-t_j)\mu_j^2 \hspace{.2cm}, 
&    d_3(\delta^3) & =  (1-h)\delta^2 - \sum_{j=0}^3\pi_j^2 \hspace{.2cm}, \\
d_4(\sigma_0^4)     & =  N\cdot (\delta^3 - \sum_{j=0}^3r_{j-1}\sigma_j^3) \hspace{.2cm}.
&   & 
\end{array}$$
We define $C_n \approx C_{n-4}$ for $n\geq 5$ with appropriate subscripts.
\end{Pro}

It is instructive to compare this resolution with the case $|G| = \infty,$ treated in
\mbox{\cite{bhzz1,bhzz2,bz}}, for which $C_j = 0$,  $j\geq 4$. Here the finiteness
of $G$ is reflected by the new class $\sigma_0^4 \in C_4$ whose boundary generates
Ker$(d_3)$, which is no longer $\{ 0\}$. The diagonal $\Delta$, taken from these same
references, suffices through dimension $3$, and therefore for computations of the cup
products into dimensions $\leq 3$. Thus, the following theorems do not give the cup products   
into dimensions $\geq 4$, these will have to wait until $\Delta_4, \Delta_5, ...$ are computed
(which at present seems very difficult), or some other method applied.

\begin{Thm} \label{ringp'}
\begin{align*}
H^l(P'_{8\cdot 3^k}; {\mathbb Z}) 
 & \iso  \left\{\begin{array}{lcll}
{\mathbb Z}       & = &
 \Big\langle 1 := \left[ \sum_{j=0}^3 {\hat \sigma}_j^0 \right] \Big\rangle 
 & ,\mbox{ if } l = 0, \\
0 &  &  & ,\mbox{ if } l = 1, \\
{\mathbb Z}_{3^k}
 \hspace{.2cm}  & = &    \<\gamma_2 := [\hat{\mu}_3^2] \>              & ,\mbox{ if } l = 2, \\
0 &  &  & ,\mbox{ if } l = 3, \\
{\mathbb Z}_{8\cdot 3^k} & =   &   \<\alpha_4 :=  [\hat{\sigma}_0^4] \>   & ,\mbox{ if } l =4. 
\end{array}\right.\\
 H^l(P'_{8\cdot 3^k}; {\mathbb Z}_3) 
 & \iso  \left\{\begin{array}{lcll}
{\mathbb Z}_3       & = &
 \Big\langle 1 := \left[ \sum_{j=0}^3 {\hat \sigma}_j^0 \right] \Big\rangle 
 & ,\mbox{ if } l = 0, \\
{\mathbb Z}_3       & = & 
 \<\beta_1 := [{\hat \rho}_3^1 - {\hat \rho}_2^1] \> 
 & ,\mbox{ if } l = 1, \\
{\mathbb Z}_3       & = &
 \<\gamma_2 := [{\hat \sigma}_2^2] \> 
 & ,\mbox{ if } l = 2, \\
{\mathbb Z}_3       & = &
 \<\delta_3 := [{\hat \delta}^3]
           = -[{\hat \sigma}_0^3] = \cdots 
           = -[{\hat \sigma}_3^3] \> 
 & ,\mbox{ if } l = 3, \\
{\mathbb Z}_3 & =  &  \<\alpha_4 :=  [\hat{\sigma}_0^4] \>  & ,\mbox{ if } l = 4. 
\end{array}\right.
\end{align*}
Furthermore,  $\beta_1^2 = 0$, $\beta_1 \gamma_1 =  -\delta_3$.
\end{Thm}

\begin{Thm}
Let  $M = S^3/P'_{8\cdot 3^k}$.
\begin{eqnarray*}
H^*(M; {\mathbb Z}) & \iso & {\mathbb Z}[\beta_2,\delta_3]^{\star}/(3^k\beta_2=0) \hspace{.2cm}. 
\\
H^*(M; {\mathbb Z}_3) & \iso & {\mathbb Z}_3[\beta_1, \gamma_2, \delta_3]^{\star}
/(\beta_1^2 = 0, \beta_1\gamma_2 = -\delta_3) \hspace{.2cm}.
\end{eqnarray*}
If $p \not= 3$, then
\begin{eqnarray*}
H^*(M; {\mathbb Z}_p) & \iso & {\mathbb Z}_p[\delta_3]^{\star} \hspace{.2cm}.
\end{eqnarray*}
\end{Thm}

\subsection[The groups B(2^k(2n+1))]{The groups $B_{2^k(2n+1)}$} \label{bgroups}

The groups $B_{2^k(2n+1)}$, $k \geq 2$, $n \geq 0$, have the presentation
$$
B_{2^k(2n+1)} = \< x, y \,|\, x^{2^k} = y^{2n+1} = 1, xyx^{-1} = y^{-1} \>\ .
$$
They also have the semidirect product structure 
$B_{2^k(2n+1)} \iso C_{2n+1} \rtimes C_{2^k}$, as seen from the split 
short exact sequence
\begin{center}\begin{picture}(180, 10)
\thinlines
\put( -5,0){\makebox(0,0)[t]{\smash{$1$}}}
\put(  0,3){\vector(1,0){10}}
\put( 27,0){\makebox(0,0)[t]{\smash{$C_{2n+1}$}}}
\put( 50,-2){\makebox(0,0)[t]{\smash{\Large{$\hookrightarrow$}}}}
\put( 48.5,6){\makebox(0,0)[t]{\smash{\small{$\vartriangleleft$}}}}
\put( 80,0){\makebox(0,0)[t]{\smash{$B_{2^k(2n+1)}$}}}
\put(113,0){\makebox(0,0)[t]{\smash{\Large{$\twoheadrightarrow$}}}}
\put(113,9){\makebox(0,0)[t]{\smash{\small{$p$}}}}
\qbezier(107,1)(114,-2)(121,1)\put(104,2){\vector(-3,1){1}}
\put(113,-8){\makebox(0,0)[t]{\smash{\small{$s$}}}}
\put(130,0){\makebox(0,0)[t]{\smash{$C_{2^k}$}}}
\put(140,3){\vector(1,0){10}}
\put(155,0){\makebox(0,0)[t]{\smash{$1$}}}
\end{picture}\end{center}
where $C_{2n+1}$ is generated by $y$, $C_{2^k}$ by $x$, $p(y) = 1$, $p(x) = x$,
$s(x) = x$.  Furthermore, $\left(B_{2^k(2n+1)}\right)_{ab} = \ZZ_{2^k}$.

As mentioned before \fullref{tblSMfundgps}, the cases  $n = 0$  or $k = 2$ reduce to groups 
that have already been
studied  (respectively  $\ZZ_{2^k}$  or $Q_{4(2n+1)}$), so we assume henceforth that
 $n \geq 1$  and $k \geq 3$. \fullref{tblSMfundgps} then gives two cases, the second and sixth, which
give $B_{2^k(2n+1)}$ as the fundamental group $\pi_1(M)$. Specifically, both
\begin{align*} M &= (O,o,0 : e ; (2,1),(2,1),(a_3,b_3)), \ {\rm with} \ a_3 = 2n + 1, \ 2^{k-2}{=} (e+1)a_3 + b_3 \ ,
\\ M' &= (O,n,1 : e ; (a_1,b_1)), \  {\rm with} \ \ a_1 = 2^{k-2}, \ \ 2n + 1 = e a_1 + b_1 \ , \end{align*}
have  $B_{2^k(2n+1)} $  as fundamental group for  $k\geq 3$. We will choose $M$ for the 
computations in this subsection, and briefly remark about $M'$ in \fullref{rem5.7} below.

Indeed, choosing $M$, the resulting resolution is formally identical to that in \fullref{prop5.2} 
(with different structure constants $e, a_j, b_j$). And the diagonal $\Delta$ is similarly taken
from \cite{bhzz1,bhzz2,bz}. The results are as follows.
 
\begin{Thm} \label{ringb}
\begin{align*}
 H^l(B_{2^k(2n+1)}; {\mathbb Z}) 
 & \iso  \left\{\begin{array}{lcll}
{\mathbb Z}       & = &
 \Big\langle 1 := \left[ \sum_{j=0}^3 {\hat \sigma}_j^0 \right] \Big\rangle 
 & ,\mbox{ if } l = 0, \\
0 &  &  & ,\mbox{ if } l = 1, \\
{\mathbb Z}_{2^k}
 \hspace{.2cm}  & = &    \<\gamma_2 := [\hat{\mu}_2^2] \>              & ,\mbox{ if } l = 2, \\
0 &  &  & ,\mbox{ if } l = 3, \\
{\mathbb Z}_{(2n+1)\cdot 2^k} & =   &   \<\alpha_4 :=  [\hat{\sigma}_0^4] \>   & ,\mbox{ if } l =4. 
\end{array}\right.
\\
 H^l(B_{2^k(2n+1)}; {\mathbb Z}_2) 
 & \iso  \left\{\begin{array}{lcll}
{\mathbb Z}_2       & = &
 \Big\langle 1 := \left[ \sum_{j=0}^3 {\hat \sigma}_j^0 \right] \Big\rangle
 & ,\mbox{ if } l = 0, \\
{\mathbb Z}_2       & = & 
 \<\beta_1 := [{\hat \rho}_2^1 - {\hat \rho}_1^1] \> 
 & ,\mbox{ if } l = 1, \\
{\mathbb Z}_2       & = &
 \<\gamma_2 := [{\hat \sigma}_2^2] \> 
 & ,\mbox{ if } l = 2, \\
{\mathbb Z}_2       & = &
 \<\delta_3 := [{\hat \delta}^3]
           = -[{\hat \sigma}_0^3] = \cdots 
           = -[{\hat \sigma}_3^3] \> 
 & ,\mbox{ if } l = 3, \\
{\mathbb Z}_2 & =  &  \<\alpha_4 :=  [\hat{\sigma}_0^4] \>  & ,\mbox{ if } l = 4. 
\end{array}\right.
\end{align*}
\end{Thm}

\begin{Thm}\label{thm5.6}
Let  $M = S^3/B_{ 2^k (2n+1)}$.
\begin{eqnarray*}
H^*(M; {\mathbb Z}) & \iso & {\mathbb Z}[\beta_2,\delta_3]^{\star}/(2^k\beta_2=0) \hspace{.2cm}. 
\\
H^*(M; {\mathbb Z}_2) & \iso & {\mathbb Z}_2[\beta_1, \gamma_2, \delta_3]^{\star}
/(\beta_1^2 = 0, \beta_1\gamma_2 = \delta_3) \hspace{.2cm}.
\end{eqnarray*}
\end{Thm}

\begin{Rem}\label{rem5.7} If the above calculations are done using the resolution based on 
the manifold $M'$ instead of $M$, the resultant cohomology rings are 
isomorphic to those given in \fullref{ringb} and \fullref{thm5.6} but with different 
generators for the cohomology. Namely, following the notation in 
\cite{bz}, the classes $\beta_1, \gamma_2$ are replaced respectively by
classes $\theta, \varphi$. It would be interesting to know whether $M$ and
$M'$ are homeomorphic.
\end{Rem}

\section{Applications and further questions} \label{sec:6}

In this section, we give a brief description of some applications of the
cohomology ring calculations in \fullref{sec:3}--\fullref{sec:5} to spherical
space forms.  We conclude with some questions and potentially interesting
directions for further research, including the $Q(8n, k, l)$ groups.

The first application (not new) is to existence of a degree $1$ map of
an orientable, closed, connected $3$--manifold $M$ to $\RR P^3$,
which in turn is related to the theory of relativity.  Indeed, it was shown 
in Shastri, Williams and Zvengrowski \cite{swz} that the
homotopy classes $[M, {\mathbb R}P^3]$ of maps from $M$ to the real
projective $3$--space ${\mathbb R}P^3$ are bijectively equivalent to the homotopy
classes of Lorentz metric tensors over the $4$--dimensional space-time
manifold $M \times {\mathbb R}$.

Let $M$ be a closed orientable connected $3$--manifold.  We say that
$M$ is of \textit{type 1} if it admits a degree $1$ map onto real
projective $3$--space ${\mathbb R}P^3$; otherwise, it is of \textit{type 2}.

We have the following theorem from \cite{swz}:
\begin{Thm}
Let $M$ be a closed orientable connected $3$--manifold.  The following are
equivalent:\begin{enumerate}
\item
The short exact sequence
$0 \to [M, S^3] \to [M, {\mathbb R}P^3] \to H^1(M; {\mathbb Z}_2) \to 0$
does not split.
\item
The manifold $M$ is of type $1$.
\item
There exists $\alpha \in H^1(M; {\mathbb Z}_2)$ with $\alpha^3 \neq 0$.
\end{enumerate}
\end{Thm}

Among the Clifford--Klein space forms, the paper \cite{sz} determined
all those having type~$1$; ie admitting a degree $1$ map onto
${\mathbb R}P^3$.  They are the space forms corresponding to the groups
$C_{2(2m+1)}$, $Q_{16n}$, and $P_{48}$.

Of course, $\RR P^3$ can be thought of as the lens space $L(2,1)$.
Let us now consider degree $1$ maps onto $L(p, q)$, $p > 2$.  We will use
the following theorem of Hayat-Legrand, Wang and Zieschang \cite{hwz1}.

\begin{Thm} \label{app0}
Let $M$ be a closed connected orientable $3$--manifold.  Assume that there
is an element $\alpha \in H_1(M)$ of order $p > 1$ such that
the linking number $a\odot a$ is equal to
$[r/p]\in{\mathbb Q}/{\mathbb Z}$ where $r$ is prime to
$p$. Then there exists a degree-one map $f\co M \to L(p, s)$ where $s$ is
the inverse of $r$ modulo $p$.
\end{Thm}

We remark that the application of this theorem uses the mod $p$ Bockstein
homomorphism $B\co H^1(X; {\mathbb Z}_p) \to H^2(X; {\mathbb Z}_p)$, which
can be simply described as arising from the connecting homomorphism of the
long exact sequence induced by the short exact sequence
$0 \to {\mathbb Z}_p \to {\mathbb Z}_{p^2} \to {\mathbb Z}_p \to 0$
of coefficients.

\begin{Thm} \label{app1}
The spherical space form $M = S^3/P_{24}$ admits a degree one map onto 
$L(3, 1)$.
\end{Thm}
\Proof
The $1$--dimensional cohomology class $x \in H^1(M; {\mathbb Z}_3)$
is represented by the ${\mathbb Z}_3$
cocycle ${\hat b} - {\hat b'}$.  This lifts to the ${\mathbb Z}_9$
cochain denoted also by ${\hat b} - {\hat b'}$.  Now
$\delta^2({\hat b} - {\hat b'})$
$ = (-{\hat c}+2{\hat c'}) - (2{\hat c}-{\hat c'})$
$ = -3{\hat c} + 3{\hat c'}$
$ = 6{\hat c'}$, since $-{\hat c} = {\hat c'}$,
and dividing this by $3$ we obtain $B(x) = 2[{\hat c'}] = 2y$.
Thus, $x\cup B(x) = x \cup 2y = 2(x\cup y) = 2\cdot 2z = z\ ({\rm mod}\ 3)$,
and applying \fullref{app0}, completes the proof.
\end{proof}

We remark that this is related to Theorem $1.1$ of \cite{bz}.
Similarly, we can show that for $M = S^3/P'_{8\cdot3^k}$,
$B_3(x) = 0$, $k \geq 2$, hence there does not exist any degree one map 
$M\to L(3, q)$.

We now give an application to Lusternik--Schnirelmann category.  To be
clear, we speak of the normalized Lusternik--Schnirelmann category 
$\cat(X)$ of a connected topological space $X$, defined to be the 
smallest integer $n$ such that $n+1$ open sets, each contractible in $X$, 
cover $X$.  It is well known that the cup length of $X$ 
(with any coefficients) furnishes a lower bound for $\cat(X)$, while 
the dimension $n$, for $X$ a finite connected CW-complex of dimension $n$,
furnishes an upper bound.  As a simple consequence, we have the next theorem.

\begin{Thm} \label{app2}
Suppose $M$ is a type $1$ closed orientable connected $3$--manifold. Then
$\cat(M) = 3$.
\end{Thm}
\Proof
Since $M$ is of type $1$, there exists $\alpha \in H^1(M; {\mathbb Z}_2)$
with $\alpha^3 \neq 0$ which implies that $\cat(M) \geq 3$.
Since $M$ is a closed $3$--manifold, there exists a finite $3$--dimensional
CW-decomposition of $M$ which implies that $\cat(M) \leq 3$.
Combining these results, we have $\cat(M) = 3$.
\end{proof}

\begin{Cor} \label{corcat}
For $G = C_{2(2m+1)}$, $Q_{16n}$, or $P_{48}$, $\cat(S^3/G) = 3$.
\end{Cor}

\begin{Rem}
Similar results about the category of orientable Seifert manifolds
with infinite fundamental groups were obtained in \cite{bz}.
However, those results also follow from work of Eilenberg--Ganea \cite{eg}
because these manifolds are aspherical.  The present results in 
\fullref{corcat} would seem to be entirely new.
\end{Rem}

To conclude, one obvious direction for further research is to complete
all calculations for the groups $P_{120}$, $P'_{8\cdot3^k}$, and
$B_{(2n+1)2^k}$, as was done for the other finite fundamental groups.
Another interesting direction is to study the cohomology rings with
other (in particular twisted, ie nontrivial $R$--module structure)
coefficients.  It seems likely that such a study could lead to further
information about degree $1$ maps and Lusternik--Schnirelmann category,
similar to \fullref{app1} and \fullref{app2} above.  

Finally, the cohomology 
ring of the $4$--periodic groups $Q(8n, k, l)$, mentioned in \fullref{sec:2}, 
is of interest. These groups have the presentation
$$Q(8n,k,l)=\langle \ x,y,z \ | \  x^2 = y^{2n} = (xy)^2, z^{kl}=1, xzx^{-1}=z^r, 
yzy^{-1}=z^{-1} \ \rangle \  ,$$
where
$n$, $k$, $l$ are odd integers that are pairwise relatively prime,
$n > k > l > 1$, and $r \equiv -1 ({\rm mod}\ k)$,
$r \equiv 1 ({\rm mod}\ l)$.
Indeed, it is known that it suffices to consider the subfamily
$Q(8p, q) := Q(8p, q, 1)$ with $r = -1$, $p$, $q$ distinct odd primes.
For these groups, an interesting balanced presentation is given by
B\,Neumann \cite{neumann}:  
$$ Q(8p,q)= \langle \ A,B \ | \ (AB)^2 = A^{2p}, \ B^{-q}AB^q = A^{-1} \   \rangle \ .$$
A proof that the two presentations for $Q(8p,q)$ give isomorphic groups 
is given in Tomoda \cite{tomoda}.
The authors attempted, but did not succeed,
to construct a $4$--periodic resolution for $Q(8p, q)$ using the
Neumann presentation.  Of course, a demonstration that no $4$--periodic
resolution exists would give an algebraic proof that $Q(8p, q)$ cannot act freely on $S^3$
(again, as mentioned in \fullref{sec:2}, a geometric proof of this result is contained in the 
work of Perelman \cite{perelman}  and his successors).

\bibliographystyle{gtart}
\bibliography{link}

\end{document}